\documentclass[11pt]{amsart}
\usepackage[margin=1in]{geometry}
\usepackage{lipsum, bm,cleveref}
\usepackage{amsfonts}
\usepackage{graphicx}
\usepackage{mathtools}

\usepackage{epstopdf}
\usepackage[percent]{overpic}
\usepackage{algorithm}
\usepackage{algpseudocode}

\usepackage{amsopn}

\begin{document}
\title{Artificial Neural Network Solver for Time-Dependent
  Fokker-Planck Equations}
\author{Yao Li}
\address{Yao Li: Department of Mathematics and Statistics, University
  of Massachusetts Amherst, Amherst, MA, 01002, USA}
\email{yaoli@math.umass.edu}
\author{Caleb Meredith}
\address{Caleb Meredith: Department of Mathematics and Statistics, University
  of Massachusetts Amherst, Amherst, MA, 01002, USA}
\email{cmeredith@umass.edu}
\thanks{Authors are listed in alphabetical order. Yao Li is partially
  supported by NSF DMS-1813246 and DMS-2108628.}

\begin{abstract}
Stochastic differential equations play an important role in various applications when modeling systems that have either random perturbations or chaotic dynamics at faster time scales. The time evolution of the probability distribution of a stochastic differential equation is described by the Fokker-Planck equation, which is a second order parabolic partial differential equation. Previous work combined artificial neural network and Monte Carlo data to solve stationary Fokker-Planck equations. This paper extends this approach to time dependent Fokker-Planck equations. The focus is on the investigation of algorithms for training a neural network that has multi-scale loss functions. Additionally, a new approach for collocation point sampling is proposed.  A few 1D and 2D numerical examples are demonstrated. 

\end{abstract}

\maketitle

\section{Introduction}
The Fokker-Planck equation plays an important role in various applications because it describes the time evolution of a stochastic differential equation, which is widely used to study noise perturbed systems or models. Since most Fokker-Planck equations have no explicit solution, numerical Fokker-Planck solvers are necessary. Before, the main difficulty of solving a Fokker-Planck equation was that the long term stability of the Fokker-Planck solution comes from the drift term of the stochastic differential equation rather than its own boundary condition. The lack of a suitable boundary condition on the numerical domain plus the high dimensionality makes many traditional methods less effective. This problem is partially solved by the first author's series of papers \cite{FPE_2019, FPE_2022}, in which a data-driven Fokker-Planck solver is developed. An artificial neural network version of the data-driven Fokker-Planck solver for stationary Fokker-Planck equation is proposed and studied in \cite{NN_FPE}. In this paper, we will both extend the work in \cite{NN_FPE} to time dependent Fokker-Planck equations and further investigate the training methods of artificial neural networks for the neural network Fokker-Planck solver. Many unaddressed problems about neural network training and training point sampling in \cite{NN_FPE} are studied in this paper. 

The main idea of the data-driven solver is that the Fokker-Planck equation has a probabilistic representation. Hence its solution can be approximated by a Monte Carlo simulation. The data-driven solver only requires a rough Monte Carlo simulation data, which is highly noisy but can be obtained at low computational cost. One important observation is that the error term in the Monte Carlo simulation data is largely spatially uncorrelated. Therefore, the Monte Carlo simulation data can be used to guide either a classical PDE solver or an artificial neural network. The data-driven Fokker-Planck solver can be seen as a data regularization process: the noisy Monte Carlo data is regularized by the Fokker-Planck operator. The goal of training is to make the solution both fit the data and satisfy the Fokker-Planck equation (or its discretization). This idea is similar to the physics-informed neural network (PINN) \cite{Karniadakis_PINN, FPEPINN}. The main difference is that the values at collocation points are from Monte Carlo sampling rather than initial or boundary conditions. 

The loss function of the neural network Fokker-Planck solver has two parts: one comes from the Fokker-Planck operator, denoted by $L^{\mathrm{loss}}_1$, the other comes from the Monte Carlo approximation of the Fokker-Planck solution, denoted by $L^{\mathrm{loss}}_2$. Because of the low accuracy of Monte Carlo simulation, the neural loss has a multi-scale feature. In the early phase of training a randomly generated artificial neural network usually has large second order derivatives, so we have $L^{\mathrm{loss}}_1 \gg L^{\mathrm{loss}}_2$. Later in the training, $L^{\mathrm{loss}}_2$ may become the dominant term due to the error in the Monte Carlo approximation. In \cite{NN_FPE}, the problem of optimizing two loss functions at different scales was solved by an algorithm called "Alternating Adam", which alternates two Adam optimizers \cite{Adam} for two loss functions. Interestingly, later we find that this approach does not work well when the differential operator in the loss term is not fully elliptic, which includes time dependent Fokker-Planck equations and the stationary Fokker-Planck equation with degenerate elliptic term. 

Therefore, in this paper, we use time dependent Fokker-Planck equations as an example to carefully examine the methods of optimizing artificial neural networks and sampling training points. Several different methods are tested and compared. In the end we concluded that the most robust training method is the "Gradient-Based Momentum Weight" method, which gradually changes the relative weight of two loss functions based on the gradients of the two loss functions from the previous epoch. We have also tested different methods of sampling training points. In addition to sampling training points proportional to the probability density as discussed in \cite{NN_FPE}, we find that it is beneficial to concentrate collocation points (meaning training points with approximate probability density) at the initial distribution and a few selected time "anchors". We believe this is because the temporal variable is only regularized by the first order derivative, so the neural network can learn the "shape" of solution easily but needs more data at "anchors" to learn the correct scale of the solution. 

We remark that this paper is {\it not} a trivial generalization of \cite{NN_FPE}. It carefully investigates training methods of the artificial neural network with multiple loss functions at different scales. It is known that PINNs have similar issues when the data at collocation points comes from experiments \cite{Noisy_PINN}. Many non-PDE neural network trainings also need to balance training loss functions at different scales. In the examples that we have tested, the new training methods developed in this paper have superior performance to both the "Alternative Adam" proposed in \cite{NN_FPE} and the idea of trainable weight proposed in \cite{Trainable_Weight, Trainable_Weight_2, Trainable_Weight_Theory}. We expect these new discoveries to be applied to other applications in the future. 

The organization of this paper is as follows. Section 2 reviews stochastic differential equations, the Fokker-Planck equation, and the data-driven Fokker-Planck solver with both the discretization version and the neural network version. The neural network solver is described in Section 3. Section 4 investigates a few different ideas of training a neural network with multi-scale loss functions, which is one of the main focus of the present paper. The numerical examples are demonstrated in Section 5. Section 6 uses some numerical examples to demonstrate the improved training result from a better sampling method. Section 7 is the conclusion. An appendix discusses the implementation and hyper-parameter selection of training methods in full details, as well as the training point  selection. 

\section{Preliminary}
\subsection{Stochastic Differential Equations and the Fokker-Planck equation}
Consider a stochastic differential equation (SDE) that is of the form 
\begin{equation}
\label{SDE}
    d\boldsymbol{X}_t = f(\boldsymbol{X}_t)dt + \sigma(\boldsymbol{X}_t)\mathrm{d}\boldsymbol{W}_t \,,
\end{equation}
where $f$ is a vector field in $\mathbb{R}^d$, $\sigma$ is a $d\times m$ matrix-valued function, and $\boldsymbol{W}_t $ is the standard Wiener process in $\mathbb{R}^m$. The solution of \eqref{SDE}, denoted by ${\bm X} = \{ {\bm X}_t \,|\, t \in \mathbb{R} \}$, is a stochastic process on $\mathbb{R}^d$. Since the theme of this paper is about the numerical method, throughout this paper, we assume that $f$ and $\sigma$ have sufficient regularities such that equation \eqref{SDE} admits a weak solution. It is well known that ${\bm X}_t$ is a continuous-time Markov process with an infinitesimal generator $\mathcal{L}$ satisfying

\begin{equation}
    \mathcal{L} h = -\sum_{i = 1}^n f_i h_{x_i} + \frac{1}{2}\sum_{i,j = 1}^{n} D_{i,j} h_{x_i, x_j} \,,
\end{equation}
where $D = \sigma^T \sigma$ is a $d\times d$ matrix-valued function.

The Fokker-Planck equation is a parabolic partial differential equation (PDE) that describes the time evolution of the probability density function of an SDE. More precisely, let $u = u(t,\boldsymbol{x})$ be the probability density function of the solution ${\bm X}_t$ to equation \eqref{SDE}, such that $u(t,\boldsymbol{x})$ is the probability density at $\boldsymbol{x} \in \mathbb{R}^d$ at time $t$. Let $D=\sigma^T\sigma$ be the diffusion matrix. The Fokker-Planck equation reads
\begin{equation}
    \label{FPE}
    u_t =  \mathcal{L}^* u = -\sum_{i=1}^{n}(f_i u)_{x_i} + \frac{1}{2} \sum_{i,j=1}^{n}(D_{i,j}u)_{x_i,x_j}  \,,
\end{equation}
where $\mathcal{L}^*$ is the adjoint operator of the generator $\mathcal{L}$. 

In addition, if the SDE \eqref{SDE} admits an invariant probability measure $\pi$, then the probability density function of $\pi$, denoted by $u( {\bm x})$, must satisfy the stationary Fokker-Planck equation, which is given by 

\begin{equation}
    \label{A}
    \mathcal{L}^*u=0 \text{ and } \int_{\mathbb{R}^d}ud\boldsymbol{x} = 1
    \,,
\end{equation}

\subsection{Data-Driven Stationary Fokker-Planck Equation Solver} The Fokker-Planck solver studied in this paper is based on the data-driven solver for stationary Fokker-Planck equations described in \cite{FPE_2019}. As discussed in the introduction, when solving the Fokker-Planck equation, many traditional PDE solvers have problems with unbounded domains and high dimensionality, while Monte Carlo simulations usually have accuracy issues. This problem is partially solved by the data-driven hybrid method proposed in \cite{FPE_2022}, which considers local Fokker-Planck equations on a subset of the entire domain without the knowledge of the boundary condition. Instead, Monte-Carlo simulation is used to provide a reference solution that makes up for the lack of a boundary condition. 

Take the 2D stationary Fokker-Planck as an example. Let $D=[a_0,b_0]\times[a_1,b_1]$ be the numerical domain, which is further split into an $N \times M$ grid of boxes. The numerical solution $\boldsymbol{u} \in \mathbb{R}^{N \times M}$ of the stationary Fokker-Planck equation is an approximation of the probability density of $u$ at the center of each grid box. Now let $\boldsymbol{A}$ represent the discretization of the operator $\mathcal{L}^*$ on $D$ with respect to all interior boxes. Then $\boldsymbol{A}$ is an $(N-2)(M-2) \times (NM)$ matrix that provides the linear constraint on $\boldsymbol{u}$ given by

\begin{equation}
    \boldsymbol{Au}=\boldsymbol{0}
\end{equation}

Next, we run a long trajectory of $\boldsymbol{X}_t$ and count the sample points in each grid box, which gives an approximated invariant probability density function denoted by $\boldsymbol{v}=\{v_{i,j}\}_{i=1,j=1}^{i=N,j=M}$
The numerical solution $\boldsymbol{u}$ is then given by the optimization problem

\begin{equation}
\begin{dcases}
    \min_{\boldsymbol{u}}||\boldsymbol{u}-\boldsymbol{v}||_2\\
    \textrm{subject to }\boldsymbol{Au}=\boldsymbol{0}
\end{dcases}
\end{equation}

It is further proved in \cite{FPE_2022} that the error in the reference solution $\boldsymbol{v}$ is significantly removed by the projection in solving the optimization problem. 

The data-driven solver for stationary Fokker-Planck equations has an artificial neural network version proposed in \cite{NN_FPE}. The idea is that the constrained optimization problem above can be replaced by a unconstrained optimization problem
\begin{equation}
    \label{unc}
    \min_{\boldsymbol{u}} ||\boldsymbol{Au}||^2_2 + ||\boldsymbol{u}-\boldsymbol{v}||^2_2
\end{equation}
that preserves the key numerical properties of the original data-driven solver in \cite{FPE_2019}. This motivates us to represent $\boldsymbol{u}$ by an artificial neural network $\tilde{u}(\boldsymbol{x}_i,\boldsymbol{\theta})$, where $\boldsymbol{\theta}$ are the trainable parameters. Since artificial neural networks are differentiable, we can further replace the discretized operator ${\bm A}$ by the differential operator $\mathcal{L}^*$. Instead of the whole domain, the optimization problem is solved with respect to a set of training points. 

Mimicking the unconstrained optimization problem in \eqref{unc}, a loss function $\bar{L}(\boldsymbol{\theta})$ is given by 
\begin{equation}
    \label{stat}
    \bar{L}(\boldsymbol{\theta})=\frac{1}{N^X}\sum_{i=1}^{N^X}(\mathcal{L}^*\tilde{u}(\boldsymbol{x}_i,\boldsymbol{\theta}))^2+\frac{1}{N^Y}\sum_{j=1}^{N^Y}(\tilde{u}(\boldsymbol{y}_j,\boldsymbol{\theta})-v(\boldsymbol{y}_j))^2
\end{equation}
where $\boldsymbol{\theta}$ represents the neural network parameters that can be updated during training, $\tilde{u}(\boldsymbol{x},\boldsymbol{\theta})$ is the neural network approximation for the probability density at $\boldsymbol{x}$ for specific parameters $\boldsymbol{\theta}$, $\boldsymbol{x}_i,\boldsymbol{y}_j \in \mathbb{R}^n$ for $i \in 1,2,...,N^X$ and $j \in 1,2,...,N^Y$ are training points without Monte Carlo approximation and collocation points with Monte Carlo approximation respectively, and  $v(\boldsymbol{y}_i)$ are the Monte Carlo approximations for the probability density at those collocation points.

\section{Neural Network solver for time-dependent Fokker-Planck Equations}

The general idea behind our artificial neural network solver for time-dependent Fokker-Planck equations largely resembles the stationary case, although the implementation and training details have many differences. 

Consider the initial value problem 
\begin{align}
    \label{IVP}
\begin{dcases}
   u_t  =  \mathcal{L}^* u = -\sum_{i=1}^{n}(f_i u)_{x_i} + \frac{1}{2} \sum_{i,j=1}^{n}(D_{i,j}u)_{x_i,x_j} \\
    u(0, {\bm x}) = u_0( {\bm x}) 
\end{dcases}
\end{align}
Similarly to the stationary case, the Fokker-Planck equation is defined on an unbounded domain and the only boundary condition is that $u(0, {\bm x})$ vanishes at infinity. Let $[0, T] \times D \subset \mathbb{R} \times \mathbb{R}^d$ be the numerical domain that we are interested in. Let $\tilde{u}(t, {\bm x}, {\bm \theta})$ be the neural network approximation of the time dependent solution to the Fokker-Planck equation \eqref{IVP}. Let $\bm{\mathfrak{X}}:=\{(t_i, \boldsymbol{x}_i) \in [0, T] \times D \mid i \in 1,2,\dots,N^X\}$ be training points without Monte Carlo approximation, $\bm{\mathfrak{Y}}:=\{(t_j, \boldsymbol{y}_j) \in [0, T] \times D \mid j \in 1,2,\dots,I^Y, I^Y+1, \dots,N^Y\}$ be collocation points from the initial distribution for $j\leq I^Y$ and collocation points with Monte Carlo approximation for $j > I^Y$, and $v(t_j, {\bm y}_j)$ for $j \in 1,2, \dots, N^Y$ be $u_0({\bm y}_j)$ if $j\leq I^Y$ and Monte Carlo approximation of the probability density at collocation point $(t_j, {\bm y}_j)$ for $j>I^Y$. Similarly to \eqref{stat}, we attempt to minimize the optimization function 
\begin{equation}
    \label{FPEloss}
    L(\boldsymbol{\theta})=\frac{1}{N^X}\sum_{i=1}^{N^X}(\mathcal{L}^*\tilde{u}(t_i, \boldsymbol{x}_i,\boldsymbol{\theta})-\tilde u_t(t_j,\boldsymbol{x}_i,\bm{\theta}))^2+\frac{1}{N^Y}\sum_{j=1}^{N^Y}(\tilde{u}(t_j, \boldsymbol{y}_j,\boldsymbol{\theta})-v(t_j, \boldsymbol{y}_j))^2 := L^{\mathrm{loss}}_1 + L^{\mathrm{loss}}_2
\end{equation}

Below we will address three key components of the neural network Fokker-Planck solver, i.e.,  the selection of collocation points, the Monte Carlo simulation that provides a reference solution, and the training of the artificial neural network. 

\subsection{Sampling collocation points}

To train the neural network, we must first sample points for $\bm{\mathfrak{X}}$ and $\bm{\mathfrak{Y}}$.  The standard method is based on the sampling method used in \cite{FPE_2022}, and can be used for sampling both $\bm{\mathfrak{X}}$ and $\bm{\mathfrak{Y}}$.  It consists of two parts, sampling uniformly across the entire numerical domain, and sampling proportional to the probability density function. Due to the fact that the density tends to concentrate near global attractors of the deterministic part of the SDE, solely uniform sampling may not be sufficient, as too many points may be chosen from low density regions. On the other hand, solely sampling according to the probability density leaves scarce points in low density regions, which can cause notable error. This can be resolved by sampling $\alpha\%$ of the points uniformly and $(1-\alpha)\%$ of the points proportional to density for some $\alpha \in [0,1]$. To facilitate Monte Carlo approximation of the probability density function, when sampling $\bm{\mathfrak{Y}}$, we move the collocation point $(t_i, {\bm y}_i)$ to the center of the grid $h$-box it is in if $i > I^Y$. The pseudocode for this algorithm can be seen below in \cref{alg:001}. For simplicity, assume that $D=[a_1,b_1]\times[a_2,b_2]\times \dots \times [a_d,b_d]$ has been split into a grid of boxes with side length $h$, and $[0,T]$ has been discretized with time steps of $\delta t$. 

\begin{algorithm}
\textbf{Input:} $\alpha \in [0,1]$, $\delta t$, $N^X$ or $N^Y$ and $I^Y$. \\
\textbf{Output:} $\bm{\mathfrak{X}}$ or $\bm{\mathfrak{Y}}$. 
\begin{algorithmic}[1]
\caption{Training and Collocation Point Sampling}\label{alg:001}
\If{$\bm{\mathfrak{X}}$}
\State{$M=N^X$, $t_0 = 0$, $i_1=1$}
\Else
\State{$M=N^Y$, $t_0 = \delta t$, $i_1=I^Y+1$}
\State{Sample $I^Y$ points from initial distribution for $\bm{V}_1$ through $\bm{V}_{I^Y}$.}
\EndIf
\For{$i=i_1$ to $i=M$}
\State Uniformly sample $c_i \in [0,1]$
\If{$c_i<\alpha$}
\State{Uniformly sample $(t_i, {\bm r}_i) \in [0, T] \times D$ with ${\bm r}_i=(r_1, \dots, r_d)$}
\If{$\bm{\mathfrak{X}}$}
\State{Set $\bm{V}_i = (t_i, {\bm r}_i)$}
\Else
\State{Set $t_i$ to closest lower multiple of $\delta t$ with $t_i= \lfloor \frac{t_i}{\delta t} \rfloor \delta t$.}
\State{Set ${\bm r}_i$ to the center of the $h$-box it belongs to with $r_j=\lfloor \frac{r_j-a_j}{h} \rfloor\*h+a_j+\frac{h}{2}$}
\State {Set $\bm{V}_i = (t_i, {\bm r}_i)$}
\EndIf
\Else
\State{Uniformly sample $t_i \in [t_0,T+\delta t]$}
\State{Set $t_i$ to closest lower multiple of $\delta t$ with $t_i= \lfloor \frac{t_i}{\delta t} \rfloor \delta t$.}
\State{Run a numerical trajectory of the SDE to time $t_i$.}
\State{Let ${\bm r}_i = \bm{X}_{t_i}=(r_1, \dots, r_d)$}
\If{$\bm{\mathfrak{X}}$}
\State{Set $\bm{V}_i = (t_i, {\bm r}_i)$}
\Else
\State{Set ${\bm r}_i$ to the center of the $h$-box it belongs to with $r_j=\lfloor \frac{r_j-a_j}{h} \rfloor\*h+a_j+\frac{h}{2}$}
\State {Set $\bm{V}_i = (t_i, {\bm r}_i)$}
\EndIf
\EndIf
\EndFor
\State{Return $\bm{V}$}
\end{algorithmic}
\end{algorithm}

An alternative of \cref{alg:001} is the "anchor method", which is introduced in this paper and explored in more detail in Section 6. The motivation of the anchor method is that the spatial and temporal variables of the solution are regularized by the first and second order derivatives respectively. The second order derivative is more sensitive against changes of network parameters. Hence the spatial variables are much easier to train. Since the solution of the Fokker-Planck equation is uniquely determined by its initial value, training $L^{\mathrm{loss}}_2$ with just points from the initial distribution and $L^{\mathrm{loss}}_1$ with the standard set of points (See Section \ref{SEC:3-3} for details) can usually determine the shape of the solution. However, the scale of the solution is usually less accurate the further away from the initial distribution it is calculated because the temporal variable is less regularized. To resolve this, one only needs a relative small number of collocation points with probability density approximations to "anchor" the scale of the solution. To update \cref{alg:001} for use with the anchor sampling method with one set of anchoring points, we simply let $t_i = T$ for $i > I^Y$ ($\delta t = T$) and set $I^Y$ close to $N^Y$ when sampling collocation points in $\bm{\mathfrak{Y}}$. The rest is identical to what is described in \cref{alg:001}.

\subsection{Monte Carlo Simulation} In order to run the neural network Fokker-Planck solver, an estimate of the probability density at the collocation points $(t_j, {\bm y}_j)$ is essential. Generally speaking the number of collocation points does not need to be very large. They only provide the role of "anchoring" the solution at the right place, while $L^{\mathrm{loss}}_1$ drives the neural network approximation to the solution of the Fokker-Planck equation. 

The first step is to sample the initial distribution $u_0( {\bm x})$, which is done here using a rejection-based method. More precisely, a random variable ${\bm y} \in \mathbb{R}^d$ is sampled uniformly from the numerical domain $D$. Next, an auxiliary random variable $p$ is uniformly sampled from $(0, \sup_{ {\bm x} \in D} u_0({\bm x}))$. The sample ${\bm y}$ is accepted if and only if $p < u_0( {\bm y})$. The process is repeated until an initial value ${\bm y}$ is accepted. We note that if the initial distribution is highly concentrated, other methods such as an MCMC sampler can also be used. 

Next, we use a Euler-Maruyama scheme to approximate the probability densities at the collocation points away from the initial distribution. The time interval $[0, T]$ is divided into $L$ steps with $\delta t = T/L$. Let $\tau_j = j \delta t$ for $j \in 0, \dots, L$ and $\bm{X}_j = \bm{X}_{\tau_j}$. The initial value $\bm{X}_0$ is given by ${\bm y}$ sampled from the initial distribution. Then we have
\begin{equation}
{\bm{X}}_j = {\bm{X}}_{j-1} +f({\bm{X}}_{{j-1}})\delta t + \sigma({\bm{X}}_{{j-1}})\sqrt{\delta t}\mathcal{N}^m(0,1) \,,
\end{equation}
where $\mathcal{N}^m(0, 1)$ is a vector in $\mathbb{R}^m$ with each entry an i.i.d. standard normal random variable. 

The approximate probability density function can be obtained from repeatedly computing and recording numerical trajectories. As mentioned previously, all $t_j$ of the collocation points $(t_j, {\bm y}_j)$ are integer multiples of the time step $\delta t$, i.e., $t_j = j \delta t$. After simulating one numerical trajectory of the Euler-Maruyama scheme, we check whether the trajectory hits the $h$-box centering at each collocation point ${\bm y}_j$ at time $t_j$, or whether $\| \bm{X}_{j} - {\bm y}_j \|_\infty < h/2$. This process is repeated $M$ times. If we assume that the $h$-box centering at collocation point ${\bm y}_j$ is hit by realizations of $\bm{X}_t$ a total of $M_j$ times at time $t_j$, the probability density $u(t_j, {\bm y}_j)$ is approximated by $v(t_j, {\bm y}_j) = M_j h^{-d}M^{-1}$. If the dimension of $D$ is high (4 or larger), one can use kernel density estimation to improve the Monte Carlo estimation of $u(t_j, {\bm y}_j)$. But in general $M$ does not have to be very large, as the neural network solver can tolerate large spatially uncorrelated noise \cite{NN_FPE}. In practice $M = 10^6$ to $10^8$ is sufficient for most 2D and 3D problems. 

The pseudocode for the Monte Carlo simulation algorithm can be seen below in \cref{alg:002}.  The Monte Carlo probability density approximation at $(t, y_1, \dots, y_d)$ can be found at $G[\lfloor \frac{t}{\delta t}][\lfloor \frac{y_1-a_1}{h} \rfloor][\lfloor \frac{y_2-a_2}{h} \rfloor]\dots[\lfloor \frac{y_d-a_d}{h} \rfloor]$ for $t\geq \delta t$, since this algorithm does not approximate the probability density at the initial distribution where it is already explicitly known. If the dimension of $D$ is too high to use a grid $G$, the grid-free sampling method described in \cite{NN_FPE} can also be used. 

\begin{algorithm}
\textbf{Input:} $u_0(\bm{y})$, $\delta t$, $M$, $s=\sup_{\bm{x}\in D} u_0(\bm{x})$. \\
\textbf{Output:} Grid of Monte Carlo estimations $G$. 
\begin{algorithmic}[1]
\caption{Monte Carlo Probability Density Estimation}\label{alg:002}
\For{$k=1$ to $k=M$}
\State Uniformly sample ${\bm y} \in D$ and $p \in (0, s)$
\While{$u_0( {\bm y}) > p$}
\State Regenerate ${\bm y} \in D$ and $p \in (0, s)$ 
\EndWhile
\State Let $\bm{X}_0 = {\bm y}$
\For{$i=1$ to $i=L$}
\State{${\bm{X}}_i = {\bm{X}}_{i-1} +f({\bm{X}}_{{i-1}})\delta t + \sigma({\bm{X}}_{{i-1}})\sqrt{\delta t}\mathcal{N}^m(0,1)$}
\State{Denote $\bm{X}_i = (x_1, \dots, x_d)$}
\State{Compute grid position of point with $\tilde{x}_j=\lfloor \frac{x_j-a_j}{h} \rfloor$}
\If{$\bm{X}_i \in D$}
\State{$G[i][\tilde{x}_1][\tilde{x}_2]\dots[\tilde{x}_d]\mathrel{+{=}} M^{-1}h^{-d}$}
\EndIf
\EndFor
\EndFor
\State{Return $G$}
\end{algorithmic}
\end{algorithm}

\subsection{Neural Network training} \label{SEC:3-3}
The last and the most important step is to train the artificial neural network to minimize the loss function $L( {\bm \theta})$ given in \eqref{FPEloss}. Throughout this paper, we use an artificial neural network with 6 feed forward hidden layers.  The neural network architecture for the neural network has node counts given by $(d + 1)\rightarrow16\rightarrow256\rightarrow256\rightarrow256\rightarrow16\rightarrow4\rightarrow1$, where $d$ is the dimension of the phase space. Each layer uses a sigmoid activation function, and the Adam Optimizer is used for the optimization.

The loss function \eqref{FPEloss} is a combination of two loss functions 

\begin{equation}
\begin{dcases}
L^{\mathrm{\textrm{loss}}}_1(\bm{\theta}) = \frac{1}{N^X}\sum_{i=1}^{N^X}(\mathcal{L}^*\tilde u(t_j,\boldsymbol{x}_i,\bm{\theta})-\tilde u_t(t_j,\boldsymbol{x}_i,\bm{\theta}))^2\\
L^{\mathrm{loss}}_2(\bm{\theta}) = \frac{1}{N^Y}\sum_{j=1}^{N^Y} (\tilde u(t_j,\boldsymbol{y}_j, \bm{\theta}) - v(t_j,\boldsymbol{y}_j))^2\\
\end{dcases}
\end{equation}

Depending on the accuracy of Monte Carlo sampler, the loss values for these two loss functions may have very different scales. Generally speaking $L_1^{\textrm{loss}}$ is large in the beginning of the training because a randomly given neural network usually has large second order derivatives. However, $L_2^{\textrm{loss}}$ could be larger than $L_1^{\textrm{loss}}$ at the end of training if the Monte Carlo approximation $v$ is not very accurate. This property of the loss functions needs to be carefully addressed. If one simply runs the Adam optimizer for the sum $L^{\textrm{loss}}_1 + L^{\textrm{loss}}_2$, or the weighted sum $L^{\textrm{loss}}_1 + \theta L^{\textrm{loss}}_2$ for some $\theta$, then one loss function can dominate the other and yield unsatisfactory results. To resolve these issues, training algorithms needed to be developed to evenly balance the two loss functions. This is one of the main focuses of this paper, and will be addressed in detail in the next section.

\section{Training Algorithms}

This section provides an overview of a number of training algorithms introduced in this paper, as well as the motivation for their use.  Hyper-parameter selection and performance sensitivity to those values will be discussed in the Appendix.

\subsection{Alternating Adam}
The first training algorithm we will consider is Alternating Adam, the training algorithm used for the stationary case in \cite{NN_FPE}.  It will later serve as a performance benchmark.     

The idea behind this algorithm is that the Adam optimizer is scaling free. Therefore, to account for the difference in scale between the loss functions, they can be separated and alternatively trained on their own mini-batches until both of their loss values are low enough.  Doing so avoids the need to find a way to balance the two loss functions given that they are being trained in isolation.  Alternating Adam is a relatively simple algorithm to implement and served well for the stationary case where the loss function dynamics were less extreme. The pseudocode for Alternating Adam can be seen below in \cref{alg:adam}.

\begin{algorithm}
\begin{algorithmic}[1]
\caption{Alternating Adam}\label{alg:adam}
\State Initialize a neural network representation $\tilde{u}(t,\boldsymbol{x},\boldsymbol{\theta})$ with undetermined parameters $\boldsymbol{\theta}$.
\State Pick a mini-batch in $\boldsymbol{\mathfrak{X}}$, calculate the mean gradient of $L_1^{\textrm{loss}}$, and use the mean gradient to update $\boldsymbol{\theta}$.
\State Pick a mini-batch in $\boldsymbol{\mathfrak{D}}$, calculate the mean gradient of $L_2^{\textrm{loss}}$, and use the mean gradient to update $\boldsymbol{\theta}$.
\State repeat steps 2 and 3 until $L_1^{\textrm{loss}}$ and $L_2^{\textrm{loss}}$ are both small enough.
\State Return $\boldsymbol{\theta}$ and $\tilde{u}(t,\boldsymbol{x},\boldsymbol{\theta})$ for epoch with minimum $L^{\textrm{loss}}$
\end{algorithmic}
\end{algorithm}

\subsection{Fixed Weight}

The next training algorithm is the Fixed Weight algorithm, which is the simplest algorithm being introduced.  It uses a fixed weighted sum of $L_1^{\textrm{loss}}$ and $L_2^{\textrm{loss}}$ as the overall loss function, which is given by, for $\theta \in [0,1]$,

\begin{align}
    \label{Fixed}
       & L^{\textrm{loss}}(\boldsymbol{\theta})=(1-\theta) L^{\textrm{loss}}_1(\boldsymbol{\theta}) + \theta L^{\textrm{loss}}_2(\boldsymbol{\theta})
\end{align}

Using a fixed weight works well if the ratio of the gradients of two loss functions remains approximately the same, since $\theta$ can be adjusted so that the loss functions are equally influential on the overall loss.  Additionally, combining the two loss functions resolves a drawback of Alternating Adam, which is that moving along the negative gradient of one loss function can potentially result in the increase of the other loss function. However, as will be seen later, the optimal $\theta$ value is a problem specific hyper-parameter, which can be difficult to choose.

\subsection{Trainable Weight}

Due to the difficulty of selecting an optimal $\theta$ for the Fixed Weight algorithm, the remaining algorithms explore ways to adjust $\theta$ during training based on the performance from prior epochs.  If we assume that the ratio of the two loss functions remains approximately fixed, we would like our algorithms to make $\theta$ converge to the optimal $\theta$.  However, given that this assumption is rarely satisfied, we instead wish to update $\theta$ after each epoch so that it will do a better job of balancing the two loss functions for the next epoch.  


Our first attempt to do this used the loss function \eqref{Fixed}, and adjusted $\theta$ in the direction of the positive gradient $\frac{\partial L^{\textrm{loss}}}{\partial \theta}$ after each epoch. However, this approach did not look stable, as $\theta$ would simply converge to either 0 or 1, since this would allow the neural network to solely minimize one loss function at the cost of the other. In particular, since $\tilde{u}(t_j,\bm{x}_j,\bm{\theta})=0$ satisfies the Fokker-Planck equation, the neural network had a tendency to move to $\theta = 0$ and produce a zero solution. 

An alternative idea developed in \cite{Trainable_Weight_2} is to move $\theta$ following the negative gradient $\frac{\partial L^{\textrm{loss}}}{\partial \theta}$ for the loss function $L^{\textrm{loss}}(\bm{\theta})=L^{\textrm{loss}}_1(\bm{\theta})+\theta L^{\textrm{loss}}_2(\bm{\theta})$. Since $\theta$ can only increase this approach works like a constraint optimization problem: $L^{\textrm{loss}}_1( {\bm \theta}) $ is optimized under the constraint that $L^{\textrm{loss}}_2( {\bm \theta})$ is already near its optimal value. Theoretically $\theta$ should converge to a certain saddle point. We refer to \cite{Trainable_Weight_Theory} for the mathematical details. 

\subsection{Loss-Based Momentum Weight}

If $L_1^{\textrm{loss}}$ and $L_2^{\textrm{loss}}$ remain roughly proportional over all epochs, then the optimal $\theta$ for the Fixed Weight algorithm would be the loss ratio 
\begin{equation}
\label{lossratio}
\frac{L_1^{\textrm{loss}}}{L_1^{\textrm{loss}}+L_2^{\textrm{loss}}}
\end{equation}
However, in practice $L_1^{\textrm{loss}}$ and $L_2^{\textrm{loss}}$ have significant fluctuations, at the very least because of the randomly sampled mini-batches. As a result, simply updating $\theta$ to \eqref{lossratio} after each mini-batch or each epoch would seriously interrupt the training. In particular, we observed that $L_1^{\mathrm{loss}}$ drops rapidly in the beginning of training because an artificial neural network with random weights usually has very large second order derivatives. If $\theta$ is immediately updated according to the loss ratio, a rapid increase of the weight of $L_1^{\mathrm{loss}}$ may cause the training process to completely focus on optimizing $L_1^{\mathrm{loss}}$. This negative feedback loop will eventually reach the trivial solution $u(t, {\bm x}) = 0$ of the Fokker-Planck equation. Instead, we must employ some method to slow the updates of $\theta$ far enough that this feedback loop is avoided.  To do this we applied the idea of "momentum" to stabilize the change of $\theta$ during each update.     

The Loss-Based Momentum Weight algorithm can be implemented in a few different ways. An initial training period of 5 epochs of Fixed Weight training are used to stabilize the results, as well as determine a better value of $\theta$ than what was used for those 5 epochs. After that, $\theta$ is set equal to a weighted average of itself and some function of the loss ratios from the previous epochs. At the 6th epoch, the initial value $r$ could be chosen as either the average of the loss ratios from the first 5 epochs, or the loss ratio at the 5th epoch. The weight $\alpha$ is a hyper-parameter that must be determined before training. See the appendix for the discussion of suitable values of $\alpha$. \cref{alg:003} shows the pseudo code for the implementation in more detail. In practice we mainly use the 5th epoch loss ratio as the initial value of $r$ because it has less training failures. Without further specification, \cref{alg:003} takes approach (1) in line $8$.

\begin{algorithm}
\begin{algorithmic}[1]
\caption{Loss-Based Momentum Weight Method}\label{alg:003}
\State Initialize a neural network representation $\tilde{u}(t,\boldsymbol{x},\boldsymbol{\theta})$ with undetermined parameters $\boldsymbol{\theta}$.
\State Set $\theta = \theta_0$
\State Set $r = 0$
\For{epochs $= 0$ to epochs $= 4$}
\State{Train using $L^{\textrm{loss}}=(1-\theta) L_1^{\textrm{loss}} + \theta L_2^{\textrm{loss}}$}
\State{Record $L_1^{\mathrm{loss}}$ and $L_2^{\mathrm{loss}}$}
\EndFor
\State {Set $r = $ (1) the most recent $\frac{L_1^{\textrm{loss}}}{L_1^{\textrm{loss}}+L_2^{\textrm{loss}}}$ or (2) the average of $\frac{L_1^{\textrm{loss}}}{L_1^{\textrm{loss}}+L_2^{\textrm{loss}}}$ over epochs $0$ to $4$}
\For{epochs $= 5$ to epochs $= N-1$}
\State Set $\theta = \alpha \theta + (1-\alpha)r$
\State{Train using $L^{\textrm{loss}}=(1-\theta) L_1^{\textrm{loss}} + \theta L_2^{\textrm{loss}}$}
\State{Set $r=\frac{L_1^{\textrm{loss}}}{L_1^{\textrm{loss}}+L_2^{\textrm{loss}}}$}
\EndFor
\State Return $\boldsymbol{\theta}$ and $\tilde{u}(t,\boldsymbol{x},\boldsymbol{\theta})$ for epoch with minimum $L^{\textrm{loss}}$
\end{algorithmic}
\end{algorithm}

Alternatively, one can further stabilize the fluctuation of $\theta$ by taking a historic average of \eqref{lossratio} throughout all training epochs, and update $\theta$ based on this average weight ratio. Since each additional loss ratio will have a progressively smaller impact on the average loss ratio, this will make the $\theta$ updates smaller over time and further prevent the negative feedback loop. The training may take longer but the value of $\theta$ is more likely to converge. We call this the "alternative implementation" of the loss-based momentum weight method when comparing algorithms. The pseudo code for this can be seen below in \cref{alg:005}.

\begin{algorithm}
\begin{algorithmic}[1]
\caption{Alternative Implementation of the Loss-Based Momentum Weight}\label{alg:005}
\State Initialize a neural network representation $\tilde{u}(t,\boldsymbol{x},\boldsymbol{\theta})$ with undetermined parameters $\boldsymbol{\theta}$.
\State Set $\theta = \theta_0$
\State Set $r = 0$
\For{epochs $= 0$ to epochs $= 4$}
\State{Train using $L^{\textrm{loss}}=(1-\theta) L_1^{\textrm{loss}} + \theta L_2^{\textrm{loss}}$}
\State{Set $r=r+\frac{L_1^{\textrm{loss}}}{L_1^{\textrm{loss}}+L_2^{\textrm{loss}}}$}
\EndFor
\For{epochs $= 5$ to epochs $= N-1$}
\State Set $\theta = \alpha \theta + (1-\alpha)\frac{r}{\textrm{epochs}}$
\State{Train using $L^{\textrm{loss}}=(1-\theta) L_1^{\textrm{loss}} + \theta L_2^{\textrm{loss}}$}
\State{Set $r=r+\frac{L_1^{\textrm{loss}}}{L_1^{\textrm{loss}}+L_2^{\textrm{loss}}}$}
\EndFor
\State Return $\boldsymbol{\theta}$ and $\tilde{u}(t,\boldsymbol{x},\boldsymbol{\theta})$ for epoch with minimum $L^{\textrm{loss}}$
\end{algorithmic}
\end{algorithm}

\subsection{Gradient-Based Momentum Weight}

The Gradient-Based Momentum Weight algorithm is motivated by the Loss-Based Momentum Weight implementations, but does not use loss ratios for the $\theta$ updates.  Instead, it uses $\big{|}\big{|}\frac{\partial L_1^{\textrm{loss}}}{\partial \boldsymbol{\theta}}\big{|}\big{|}_2$ and $\big{|}\big{|}\frac{\partial L_2^{\textrm{loss}}}{\partial \boldsymbol{\theta}}\big{|}\big{|}_2$ in place of $L_1^{\textrm{loss}}$ and $L_2^{\textrm{loss}}$ for the Loss-Based Momentum Weight method.  This completely circumvents the feedback problem with Loss-Based Momentum Weight training, as the norms of the loss gradients do not behave the same way as the loss values themselves. In addition, since the Monte Carlo data $v(t_i, {\bm y}_i)$ has some error, the value of $L^{\mathrm{loss}}_2$ with respect to the theoretical solution is nonzero. Hence a Loss-Based Momentum Weight algorithm may cause the optimization to focus too much on $L^{\mathrm{loss}}_2$ in the late phase of training and cause over fitting. This problem is also avoided by updating the weight according to the gradient. 

Because it is very computationally expensive to compute the norms of the gradients for each mini batch, at the end of each epoch an additional batch with 500 randomly selected collocation points is run, and the gradients are computed based off of this batch. We find that averaging the ratio over the first five Fixed Weight training epochs does not make a meaningful difference. Hence throughout this paper, the initial value is chosen to be the ratio of gradients after the 5th epoch. The pseudo code for the Gradient-Based Momentum Weight algorithm can be seen below in \cref{alg:006}.  

\begin{algorithm}
\begin{algorithmic}[1]
\caption{Gradient-Based Momentum Weight}\label{alg:006}
\State Initialize a neural network representation $\tilde{u}(t,\boldsymbol{x},\boldsymbol{\theta})$ with undetermined parameters $\boldsymbol{\theta}$.
\State Set $\theta = \theta_0$
\State Set $a, b, r = 0$
\For{epochs $= 0$ to epochs $= 4$}
\State{Train using $L^{\textrm{loss}}=(1-\theta) L_1^{\textrm{loss}} + \theta L_2^{\textrm{loss}}$}
\State{Compute $a=\big{|}\big{|}\frac{\partial L_1^{\textrm{loss}}}{\partial \boldsymbol{\theta}}\big{|}\big{|}_2$ and $b=\big{|}\big{|}\frac{\partial L_2^{\textrm{loss}}}{\partial \boldsymbol{\theta}}\big{|}\big{|}_2$ based on an additional batch with 500 collocation points, and set $r=\frac{a}{a+b}$}
\EndFor
\For{epochs $= 5$ to epochs $= N-1$}
\State Set $\theta = \alpha \theta + (1-\alpha)r$
\State{Train using $L^{\textrm{loss}}=(1-\theta) L_1^{\textrm{loss}} + \theta L_2^{\textrm{loss}}$}
\State{Compute $a=\big{|}\big{|}\frac{\partial L_1^{\textrm{loss}}}{\partial \boldsymbol{\theta}}\big{|}\big{|}_2$ and $b=\big{|}\big{|}\frac{\partial L_2^{\textrm{loss}}}{\partial \boldsymbol{\theta}}\big{|}\big{|}_2$ based on an additional batch with 500 collocation points, and set $r=\frac{a}{a+b}$}
\EndFor
\State Return $\boldsymbol{\theta}$ and $\tilde{u}(t,\boldsymbol{x},\boldsymbol{\theta})$ for epoch with minimum $L^{\textrm{loss}}$
\end{algorithmic}
\end{algorithm}

\section{Numerical Example with Performance Analysis}
\subsection{1D Example}
The first numerical example studies the solution to the simple stochastic differential equation
\begin{align}
    \label{1D}
    dX_t = (-X^3_t+X_t)dt+ dW_t
\end{align}
on the numerical domain $[0,0.4] \times [-2.5,2.5]$, which is discretized into a 200 by 500 grid. Two different initial distributions are considered. The first is the standard normal distribution with probability density function
\begin{equation}
u_0(x) = \frac{1}{\sqrt{2 \pi}} e^{-\frac{1}{2}x^2} \,,
\end{equation}
the second one is a multimodal distribution given by
 \begin{equation}
     \label{multimodal}
      u_0( x) = \frac{1}{Z}(1 + \cos (5 x))e^{-\frac{1}{2}x^2} \,,
 \end{equation}
where $Z=(1+\exp{(-25/2)})\sqrt{2\pi}$ is the normalizer that makes $u_0$ a probability density function. Generally speaking, it is more difficult to train a neural network to accurately fit a multimodal distribution.

\begin{figure}[hbtp]
\centering
\label{1Dsol}
{\includegraphics[width=0.75\textwidth]{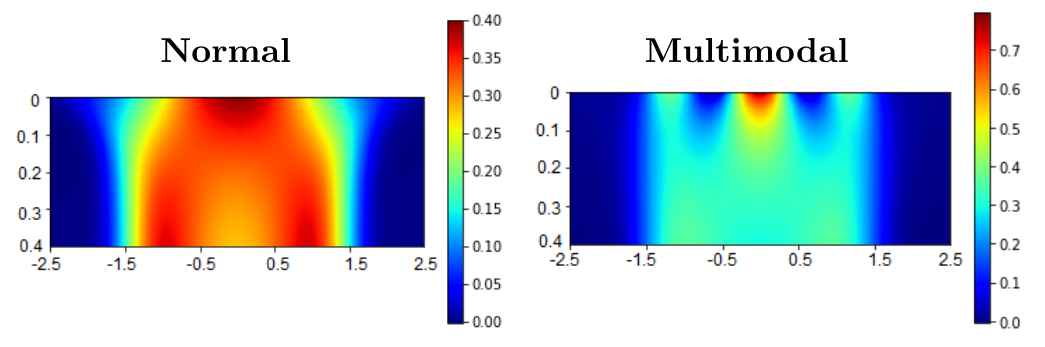}}

\caption{1D median neural network solutions from 31 samples for the normal and multimodal initial distribution using the Gradient-Based Momentum Weight method.}  
\end{figure}

The detail of algorithm implementation and parameter selection of equation \eqref{1D} is discussed in the appendix.  In summary, the size of collocation points $\bm{\mathfrak{X}}$ and training points $\bm{\mathfrak{Y}}$ are $1500$ and $20,000$ respectively. The hyper-parameters selected were $\theta=0.975$ for Fixed Weight, $\theta_0=0$ for Trainable Weight, $\theta_0=0.99$ for Loss-Based Momentum Weight with $\alpha=0.4$ for the regular implementation and $\alpha=0.6$ for the alternative implementation, and $\theta_0=0.99$ with $\alpha=0.4$ for Gradient-Based Momentum Weight.  Each of these training algorithm configurations were trained 31 times for both initial distributions, the normal distribution and the multimodal distribution, so as to demonstrate the average performance and the variation in performance for simple and complex initial distributions.  The training time is around 60-65 minutes. The neural network solutions from the Gradient-Based Momentum Weight method are shown in Figure \ref{1Dsol}.  On the left is the median solution for the normal initial distribution, and on the right is the median solution for the multimodal initial distribution.  Due to the relatively low $L^2$ error in comparison to the 2D SDE discussed later, the solutions look very similar across training methods, so only these solutions are provided as an example.

\begin{figure}[h!]
\centering
\label{1Derror}
{\includegraphics[width=0.65\textwidth]{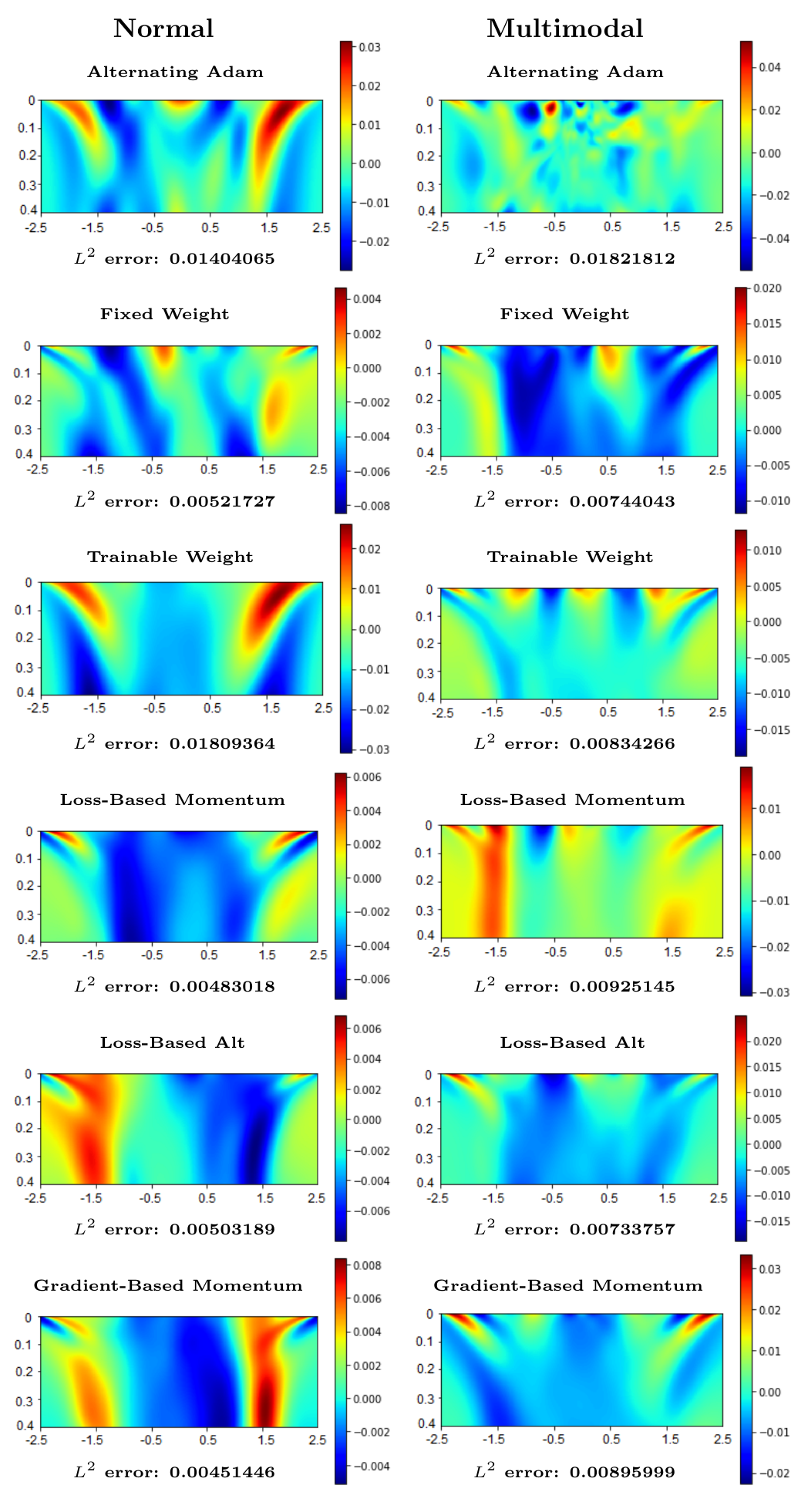}}

\caption{1D SDE median $L^2$ error heat maps for each training algorithm from 31 samples.  Left column: normal initial distribution.  Right column: multimodal initial distribution.}
\end{figure}

\subsubsection{Analysis of 1D Performance Results}
In Figure \ref{1Derror}, we demonstrate the distribution of error from all six training algorithms.  Here the ground truth solution comes from the Crank-Nicolson PDE solver on a mesh that is further refined by $4$ times. It is well known that the Crank-Nicolson scheme is a second order scheme \cite{thomas2013numerical}. Its theoretical magnitude of error on the refined mesh is around $10^{-5}$. 

\begin{figure}[h]
\label{008}
\centering
{\includegraphics[width=1\textwidth]{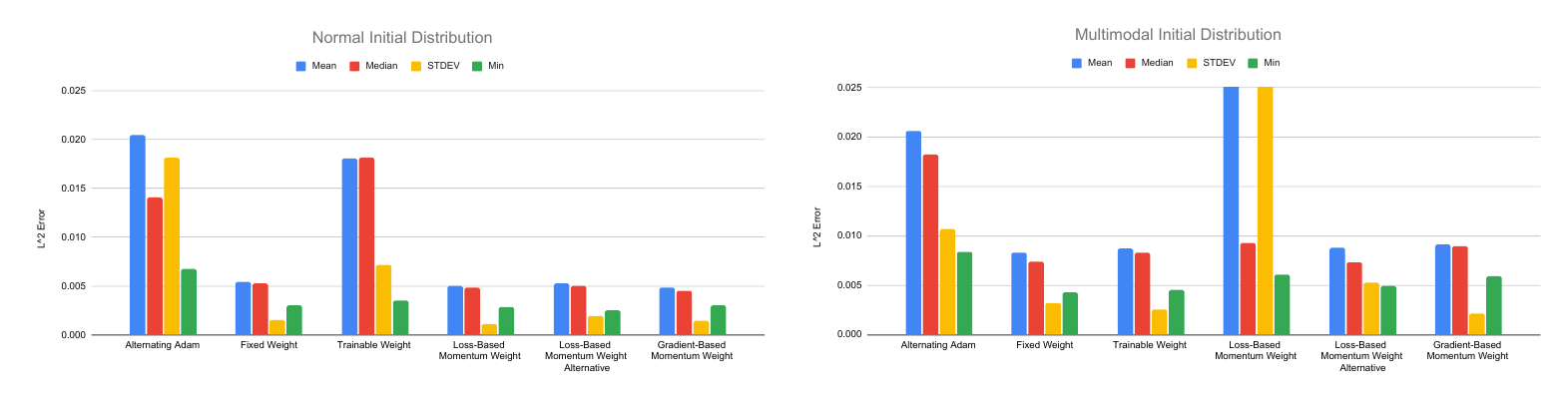}}

\caption{$L^2$ error for the normal and multimodal initial distributions. Because of the 6 failed trainings for Loss-Based Momentum Weight, the mean and standard deviation of their results extend above the graph to 0.04531 and 0.08335 respectively.}
\end{figure}

\begin{figure}[h!]
\label{011}
{\includegraphics[width=\textwidth]{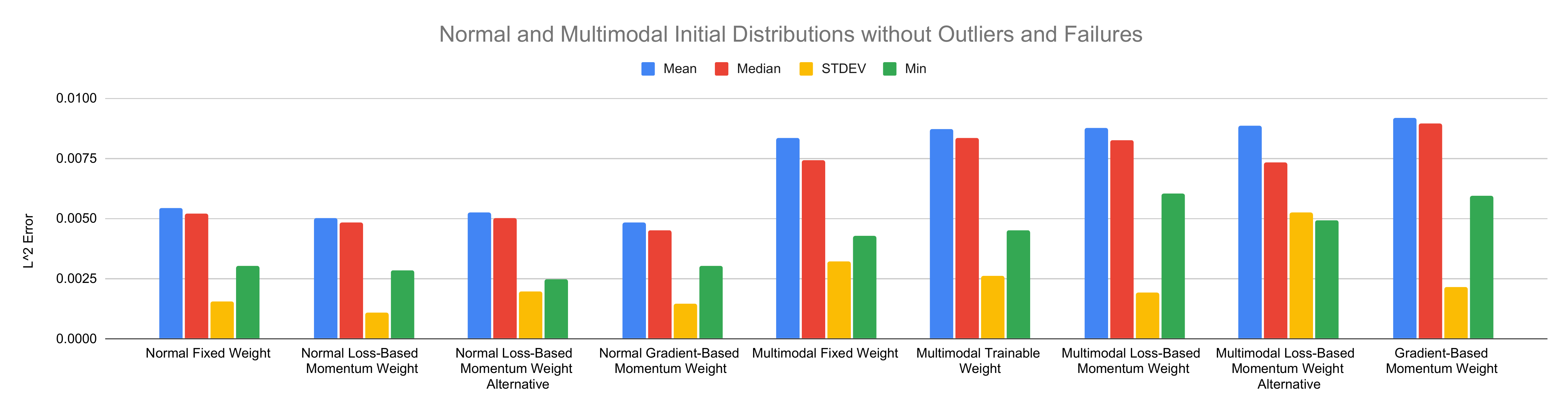}}
\caption{$L^2$ error results for both initial distributions without outliers or training failures.}
\end{figure}

The $L^2$ error of the six neural network train algorithms with the two initial distributions is demonstrated in Figure \ref{008}. It is easy to see that the performance for the simpler $\mathcal{N}(0,1)$ initial distribution is better than for the multimodal initial distribution in equation \eqref{multimodal}. This is expected because training a neural network to fit a multimodal function is more difficult. Additionally, this confirms our expectation that hyper-parameters selected for a more complicated multimodal initial distribution can be used directly on simpler initial distributions. As seen in Figure \ref{1Derror}, the error
for the multimodal initial distribution is mostly concentrated at the initial distribution, most notably near $x=-2.34$ and $x=2.34$ where the smallest peaks are typically only about half the height they should be.  While there is also error concentration there for Loss-Based Momentum Weight method for the normal initial distribution, this is far less of a problem given that the distribution should be close to zero there, and the error itself is much lower.  The more noticeable problem for the Loss-based Momentum Weight method is the maximum possible errors for the multimodal initial distribution.  This is because a portion of trainings (6 out of 31) made $\theta$ converge towards zero. Given how visually apparent these errors are, it is worth considering the performance for these methods once we discard those training results, as will be shown in Figure \ref{011}.  

Alternating Adam is a clear outlier for both initial distributions as was to be expected based on previous results.  The Trainable Weight algorithm does not work well either. Interestingly it gives higher error for the standard normal initial distribution. This likely means that the optimal Fixed Weight $\theta$ value is actually lower for the normal initial distribution because the $L_1$ error is proportionately smaller, so revisiting Fixed Weight for the normal distribution could possibly produce even better results.

For better comparison, in Figure \ref{011} we only demonstrate the error statistics of training algorithms with good results. This means that Alternative Adam for both initial distribution and Trainable Weight for the standard normal distribution are removed, as well as the 6 failed training results in the Loss-Based Trainable Weight method.

When taking both sets of results into account, combine with the hyper-parameter selection and implementation details discussed in the appendix, we can see that the Alternating Adam method is the easiest to implement but has much higher error than for the stationary Fokker-Planck equation reported in \cite{NN_FPE}. The Fixed Weight method works the best but to the best of our knowledge the weight has to be manually selected for each problem, because it is very difficult to estimate the scales of $L^{\mathrm{loss}}_1$ and $L^{\mathrm{loss}}_2$ without training the neural network. This makes the Fixed Weight method less practical. The Trainable Weight method is supported by some literature theoretically but has no performance advantage in our testing, especially considering the relatively high error for the case of the normal initial distribution. The momentum algorithms are very comparable to, or better than, the Fixed Weight method.  The Loss-Based Momentum method converges quickly but has some stability issues, meaning one must manually check whether $\theta$ converges to either $0$ or $1$. The alternative implementation has better stability, but the hyper-parameter selection process shows that it could have a slower convergence to the optimal $\theta$. (See Figure \ref{hist} in appendix.) The Gradient-Based Momentum Weight is a better balance of stability, easier implementation, and performance.  

\subsection{2D Example}
\subsubsection{2D SDE Overview and Performance Results}

Next, we apply our training methods to the same "ring example" studied in \cite{NN_FPE}. The SDE is given by
\begin{align}
    \label{ring}
    \begin{cases}
    dX_t=(-4X_t(X_t^2+Y_t^2-1)+Y_t)dt+dW_t^x \\
    dY_t=(-4Y_t(X_t^2+Y_t^2-1)-X_t)dt+dW_t^y \\
    \end{cases}
\end{align}
where $W_t^x$ and $W_t^y$ are independent one dimensional Wiener processes. It is easy to check that equation \eqref{ring} admits an explicit stationary distribution $\exp({-2(x^2+y^2-1)^2})/K$ with $K=\pi\int_{-1}^{\infty}\exp{(-2t^2)}dt$, which concentrates at the unit circle. For this example the initial distribution used is the multivariate normal distribution $\exp{(-0.5(x^2+y^2))}/(2\pi)$. Because the vector field symmetrically pushes the density to the unit circle, the solution converges to the stationary distribution quickly. We selected this equation as a numerical example because if a neural network can accurately approximate a fast-changing ring-shaped solution, we expect it can also approximate Fokker-Planck solutions in simpler shapes.

\begin{figure}[hbtp]
\label{2Dsol}
\centering
{\includegraphics[width=\textwidth]{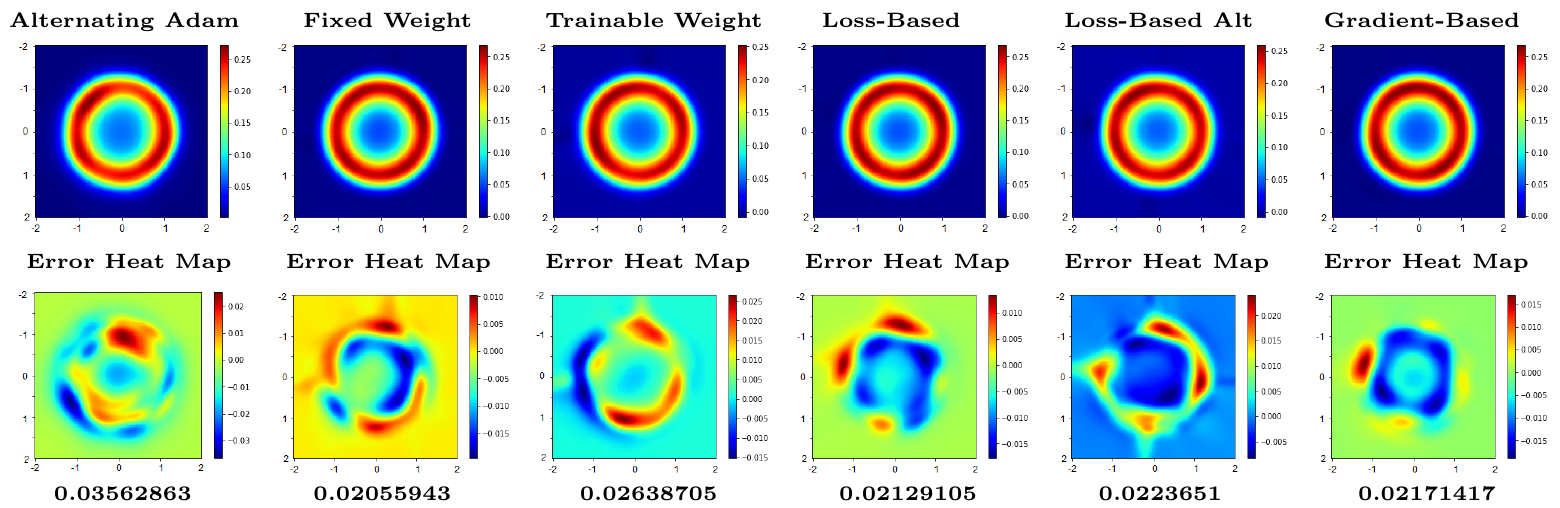}}

\caption{2D SDE median $L^2$ error results at $t=0.2$ for each training algorithm from 55 samples across training point counts.  First row: neural network output at $t=0.2$.  Second row: error heat maps, with respective $L^2$ error listed below.}
\end{figure}

The numerical domain used was $[0,0.2] \times [-2,2]\times [-2,2]$ discretized into a 200 by 200 by 200 grid.  The same neural network architecture as the 1D example was used, although the Fixed Weight training algorithm was slightly modified to have $\theta = 0.984$, which was selected based on the loss ratios from some preliminary trainings. The training data used here is also the same as the 1D case, however the average run time is around 120 minutes instead of the 60-65 minutes from before.  Each training algorithm configuration was trained 11 times at 5 different collocation point counts, namely 1083, 1875, 3468, 7500, and 13467, to see whether more points would be required for this higher dimensional SDE.  The same ratios between $N^X$, $N^Y$, and $I^Y$ were maintained. Additionally, $\bm{\mathfrak{Y}}$ is again 20,000 points sampled uniformly.  The training results are demonstrated in Figure \ref{2Dsol}. The first row shows the training result at $t = 0.2$ for the six different neural network training algorithms. The second row is the difference between neural network solution and the solution from Crank-Nicolson scheme. For each algorithm, the solution demonstrated in Figure \ref{2Dsol} is the solution with median error among the $55$ training results combined across the different training point counts, as the performance difference between them is negligible. 

Comparisons here are for the slice of the distribution when $t=0.2$, since by then the distribution is close to stationary and can be compared to the results from \cite{NN_FPE}.  Additionally, performance at $t=0.2$ is a decent indicator of overall performance, and slices at multiple times will be addressed in the next subsection.

\subsubsection{Analysis of 2D Performance Results}

\begin{figure}[h!]
\label{012}
\centering
{\includegraphics[width=1\textwidth]{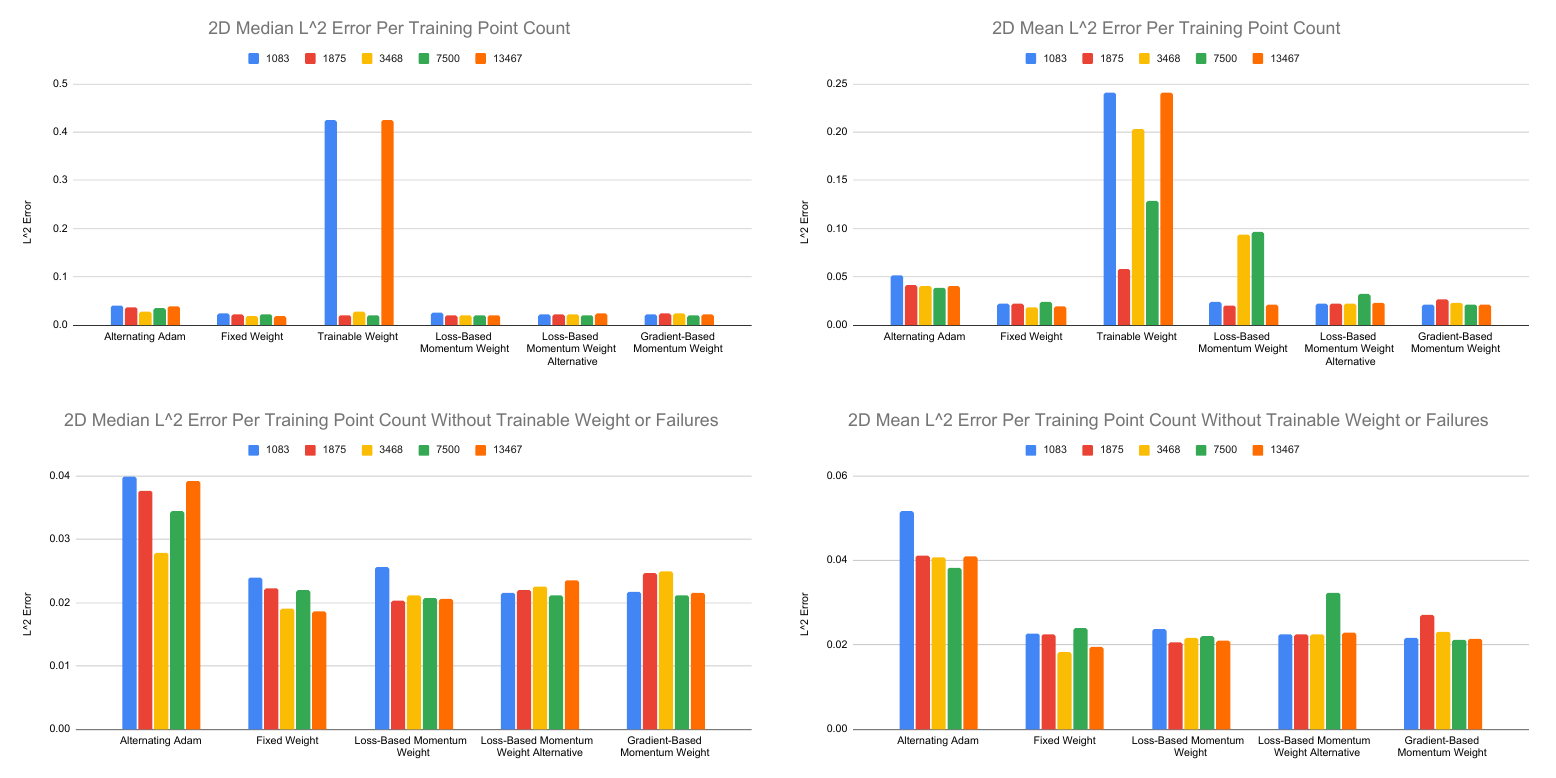}}
\caption{Median and Mean $L^2$ error for 2D example for all training algorithms using 1083, 1875, 3468, 7500 and 13467 training points.  Top: including Trainable Weight method and the Loss-Based Momentum Weight training failures.  Bottom: not including them.}
\end{figure}

Figure \ref{012} shows the median (top left) and mean (top right) $L^2$ error for each algorithm and training point count. It is easy to see that Trainable Weight is a clear outlier here and therefore should not be used. The Loss-Based Momentum Weight method also has $4$ failed trainings (out of $55$). After removing the Trainable Weight method and the failed training results from the Loss-based Momentum Weight method, a more refined result is demonstrated at the bottom of Figure \ref{012}. It is easy to see that Alternating Adam has clearly higher error, while the rest of the algorithms have rather similar performance.  

\begin{figure}[h!]
\label{013}
\centering
{\includegraphics[width=\textwidth]{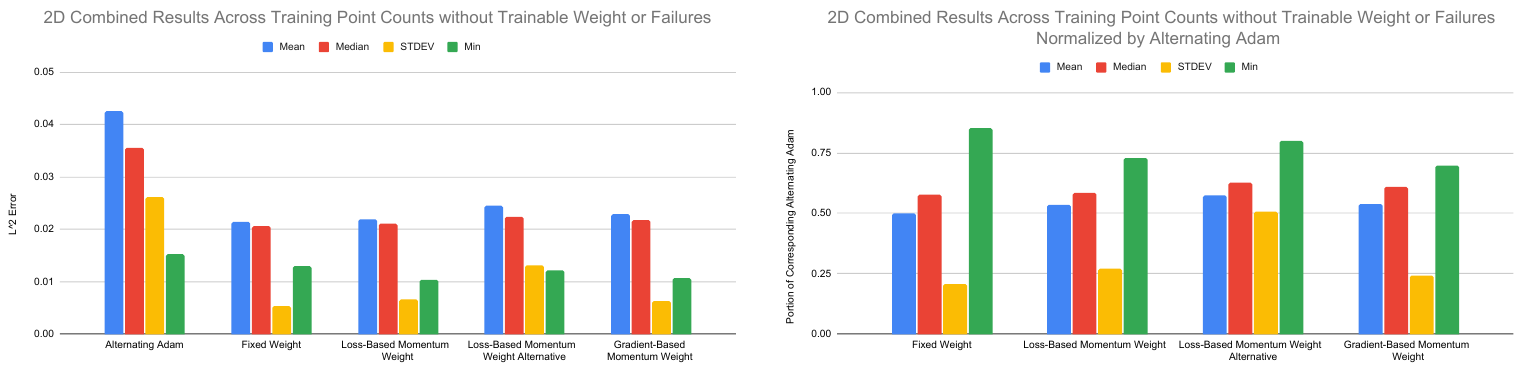}}
\caption{Median, Mean, Standard Deviation and Minimum $L^2$ error for 2D example for all training algorithms except Trainable Weight, using combined training point counts without Loss-Based Momentum Weight training failures.  Left: pure values.  Right: normalized by Alternating Adam.}
\end{figure}

Figure \ref{012} also demonstrates that there is no apparent difference between training with more than 1083 points, which is close to the maximum of 1024 training points used for this 2D ring example in \cite{NN_FPE}.  Because of this, the results are combined across training point counts in Figure \ref{013} to increase the sample size.  On the left we have the $L^2$ error values, and on the right we have those values normalized by the benchmark Alternating Adam method. Similarly to the 1D case, Fixed Weight has the lowest mean, median and standard deviation, but requires manual selection of $\theta$. The Gradient-Based Momentum Weight is better than all implementations of the Loss-Based Momentum Weight method in almost all categories. Considering all factors across the 1D and 2D cases, the Gradient-Based Momentum Weight has the best performance and will be used in our future studies, including the next section.          

\section{Comparison with Anchor Sampling Method}

\begin{figure}[hbtp]
\label{015}
\centering
{\includegraphics[width=\textwidth]{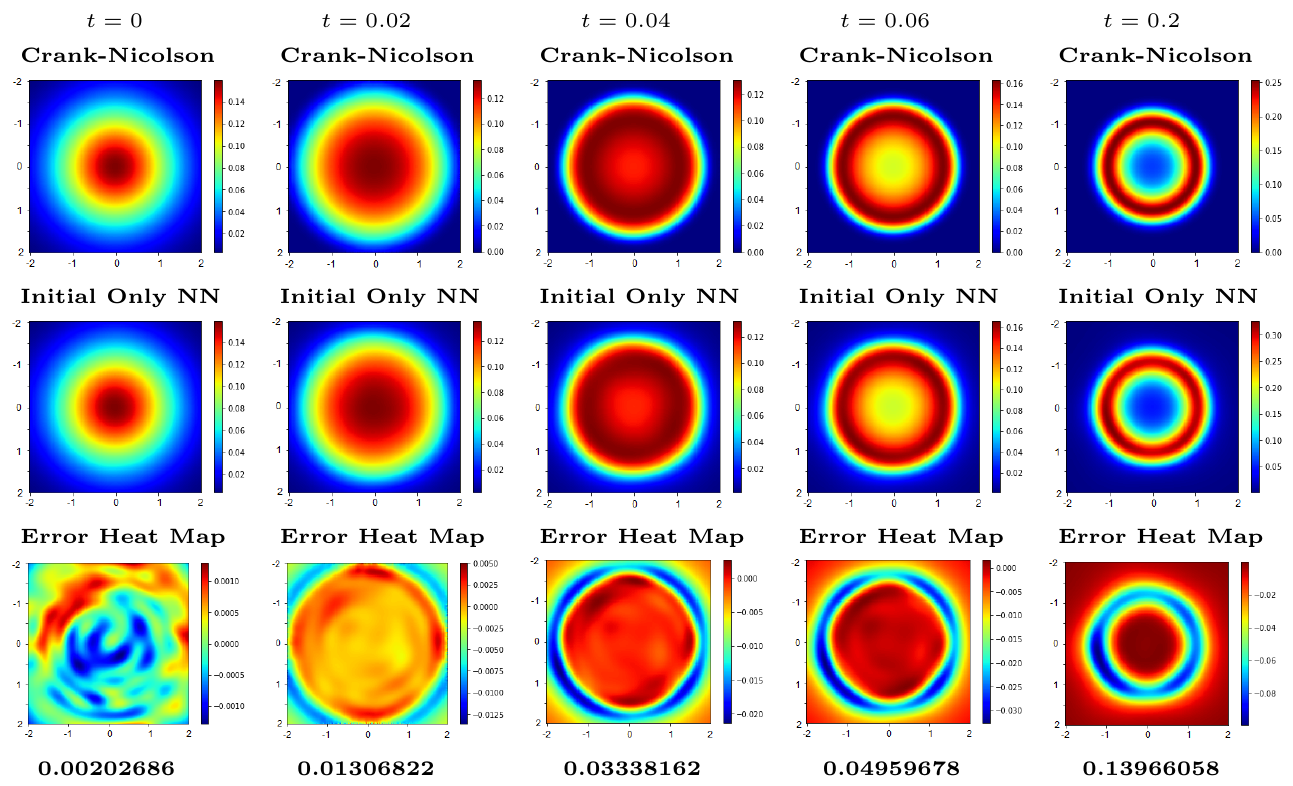}}

\caption{Crank-Nicolson solution, neural network solution and error heat maps at $t=0$, $0.02$, $0.04$, $0.06$, and $0.2$ after training the neural network with 40,000 points from the initial distribution as $\bm{\mathfrak{Y}}$.  $L^2$ error is listed below the respective heat map.}
\end{figure}

The idea behind Anchor Sampling is that theoretically the solution to the Fokker-Planck equation is uniquely determined by the initial distribution, if it is given. Therefore, we should not treat the time variable simply as "yet another dimension". Instead, we find that it is beneficial to concentrate collocation points at the initial time and the terminal time. To see this, in the 2D ring example, we let $\bm{\mathfrak{Y}}$ consist of just 40,000 points from the initial distribution and $\bm{\mathfrak{X}}$ the standard set of points uniformly sampled throughout the entire numerical domain. The result is shown in Figure \ref{015}, in which the error heat maps show the Crank-Nicolson solution minus the neural network solution. We can see that although the $L^2$ error increases with the time, the general shape of the distribution is largely preserved. Possibly due to the effect of the limited training set, the neural network does not get the scale right. Therefore, producing an accurate solution only requires correcting the scales away from the initial distribution. This motivates us to add a relatively small amount of collocation point at $t=0.2$ that serves as an "anchor".


\subsection{2D Ring SDE Anchor Sampling Numerical Results} 
\begin{figure}[h!]
\label{016}
\centering
{\includegraphics[width=1\textwidth]{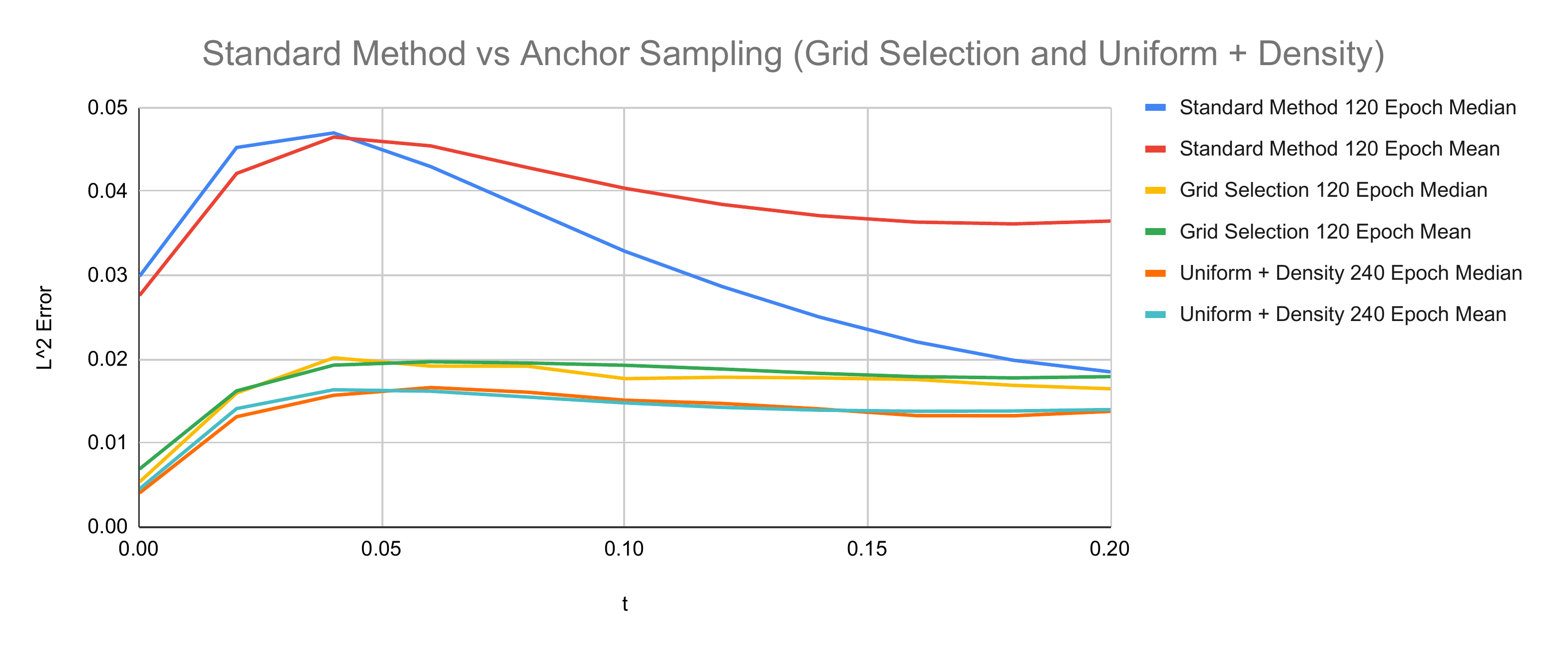}}
\caption{Median and Mean $L^2$ error comparison between Standard method for 120 epochs, Anchor sampling with Grid Selection for 120 epochs, and Anchor sampling with 0.5 Uniform Sampling 0.5 Proportional to Density Sampling for 240 epochs.  11 sample trainings used for each configuration.}
\end{figure}

\begin{figure}[hbtp]
\label{017}
{\includegraphics[width=\textwidth]{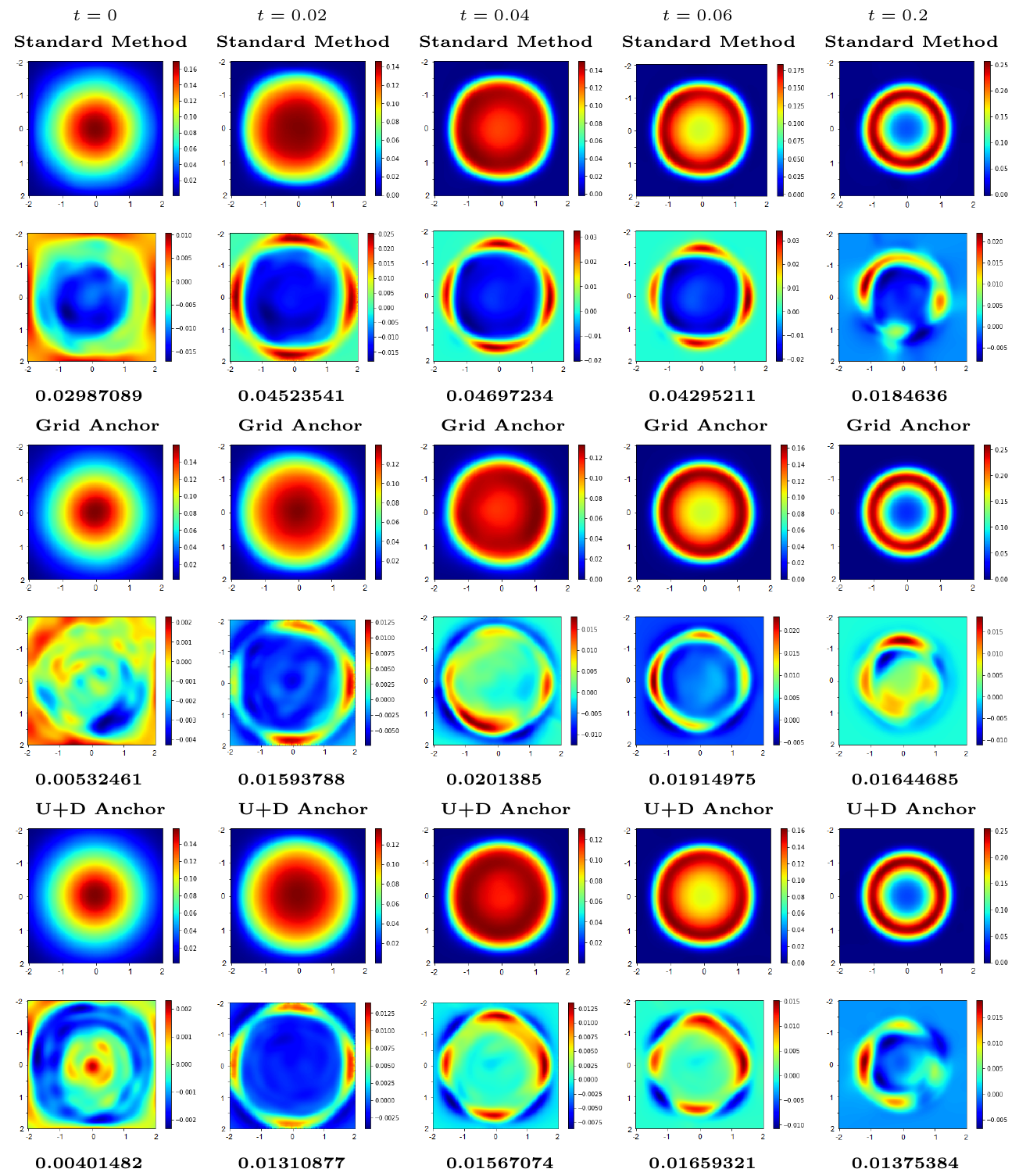}}

\caption{Neural network solutions and error heat maps of median results from 11 sample trainings per configuration at $t=0, 0.02, 0.04, 0.06, 0.2$ for the Standard method, Grid Selection Anchor sampling for 120 epochs (Grid Anchor), and Uniform + Density Anchor Sampling for 240 epochs (U+D Anchor).  $L^2$ error is listed below the respective heat map.}
\end{figure}

In the next numerical result, the set of collocation points $\bm{\mathfrak{Y}}$ is $40,000$ collocation points from the initial distribution and $1156$ points at $t = 0.2$ (see Figure \ref{anchor} in the appendix.) We tested two different ways of sampling collocation points at $t = 0.2$: one uses the standard sampling method described in Algorithm \ref{alg:001}, the other selects points from a sparser grid laid over the grid at $t = 0.2$, in this case 33 by 33 squares.  The numerical result, error heat maps, and median $L^2$ error at $t=0, 0.02, 0.04, 0.06, 0.2$ are demonstrated in Figure \ref{017}, and the median and mean $L^2$ error across all time slices is shown in Figure \ref{016}. In Figure \ref{017}, the two aforementioned sampling methods are called "U+D" and "Grid" respectively, because Algorithm \ref{alg:001} samples, in this case, half the collocation points uniformly and the other half from the probability density. The Standard Method and Grid based Anchor Sampling method were trained for 120 epochs, whereas the U+D based Anchor Sampling method was trained for 240 epochs.  This is because when the Grid based and U+D based Anchor Sampling methods were both trained at 120 and 240 epochs, one slightly outperformed the other both times.  We can see that there are diminished returns from doubling the training epochs, and that even at 120 epochs the Anchor Sampling method reduced the median $L^2$ error by about one half away from the initial distribution and the terminal time.  Due to the high concentration of points at the initial distribution for the Anchor Sampling method, the $L^2$ error there is considerably lower.  Additionally, Figure \ref{016} shows that while the median $L^2$ error for the Standard Method approaches that of the Anchor Sampling method near $t=0.2$, the mean $L^2$ error does not.  This confirms the advantage of the Anchor Sampling method. It is beneficial to use most collocation points to approximate the initial distribution, and a relative small number collocation points at the terminal time to "anchor" the solution.

\subsection{1D Multimodal SDE Anchor Sampling Numerical Results}

\begin{figure}[h!]
\label{019}
\centering
{\includegraphics[width=\textwidth]{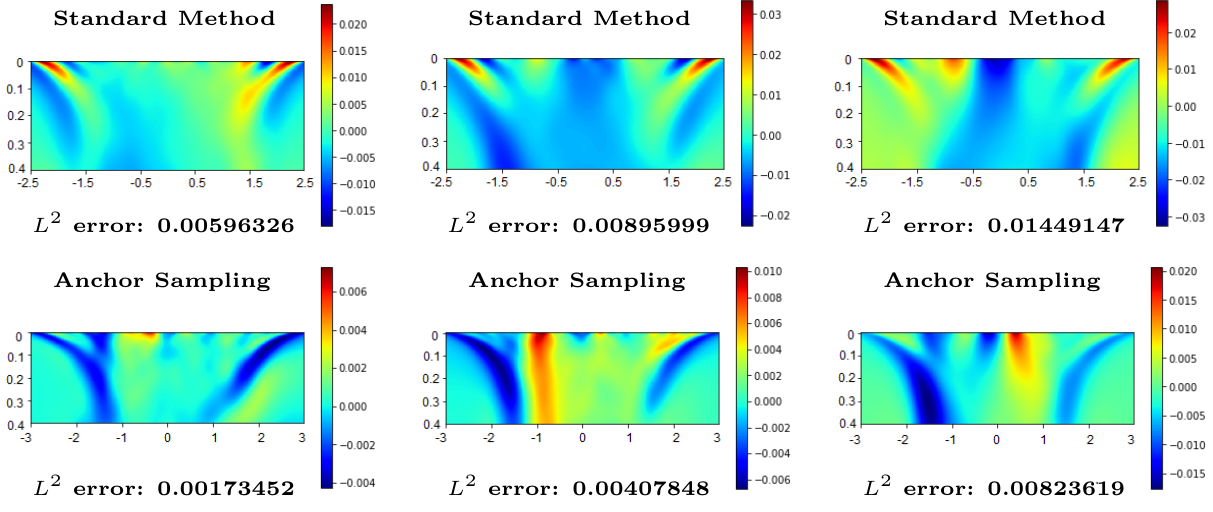}}

\caption{Minimum, Median and Maximum $L^2$ error results from 31 trainings of the 1D Multimodal SDE using the Standard Method and Anchor Sampling for 240 epochs.}
\end{figure}

For the 1D SDE studied previously, only the multimodal initial distribution was tested with Anchor Sampling, as the solution from the normal initial distribution is already satisfactory. In our computation, we sampled 2,500 points from the initial distribution and 500 points at the terminal time ($t=0.4$). Since there were only 500 grid points  at the terminal time, all of them were used and therefore no true sampling was required. The neural network was trained for 240 epochs. The result is demonstrated in Figure \ref{019}. The ground truth is still obtained from Crank-Nicolson scheme with a refined mesh. However, inspecting the range of error values, one can see that the Anchor sampling is clearly outperforming the Standard method, especially when it comes to the initial distribution. The peak of the local maximums of the initial probability density function near -2.34 and 2.34 is about 0.0432. Even in the best case scenario, the Standard method is missing about half the density there, and often misses it entirely, leaving the density at zero. In comparison, Anchor sampling is able to almost entirely eliminate that error.

\subsection{Longer Time frames}

So far we have only considered Anchor Sampling for relatively short time frames, 0.2 and 0.4 for the 2D and 1D SDEs respectively.  To extend this technique to longer time frames, multiple slices, or Anchors, are required.  This is because there is a short effective range of influence before and after each Anchor where the scale of the distributions are correct.  Additionally, the shape of the distribution starts to deform in addition to the scale drifting for sufficiently large distances from any training points.  Because of this, if just the initial distribution and terminal Anchor are used for a long time frame, error will be low near the start and end but rise considerably in between, away from the influence of either set of points. 

To demonstrate the use of multiple Anchors, the 2D ring example was trained using Anchors at $t=0.2, 0.4, 0.6, 0.8, \textrm{and } 1.0$, along with the entire initial distribution.  The $L^2$ error results from this can be seen below in Figure \ref{020}. We can see that during the time interval $[0, 0.2]$, having multiple "anchors" does not change the result very much. In addition, when multiple anchors are used, the $L^2$ error is largely consistent throughout the entire time domain.

\begin{figure}[h!]
\label{020}
{\includegraphics[width=1\textwidth]{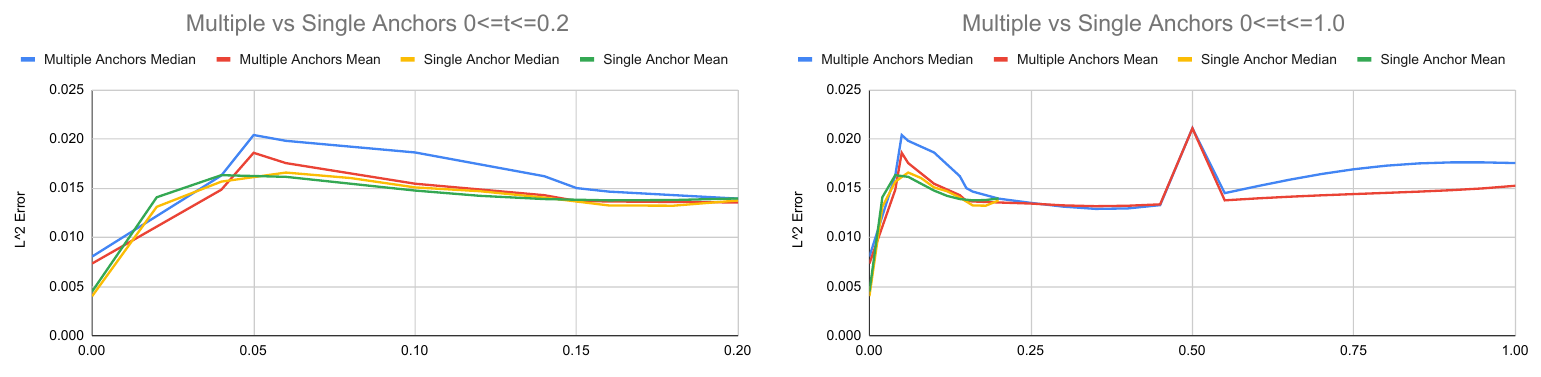}}
\caption{Comparison of mean and median $L^2$ error after 240 epochs of single vs multiple Anchors for the 2D ring SDE.}
\end{figure}

\section{Conclusion}
In this paper we examined the neural network training for the neural network Fokker-Planck solver in full detail. The main challenge here is the presence of multiple loss functions at different scales. We believe this challenge can also appear when using PINN to solve equations with noisy experimental data. One very interesting finding is that the optimization method for training the neural network seems to be problem dependent. The idea of "Alternating Adam" that worked very well for the stationary Fokker-Planck equation does not have satisfactory performance for the time-dependent Fokker-Planck equation. Instead, we tested a few different ways to balance multiple loss functions. Our analysis shows that the most robust approach is to let the relative weight of a loss function depend on the norm of the gradient of this loss function, because each update is based on the gradient rather than the value of each loss function. In addition, one needs a "momentum" term to gradually change the weight of the loss functions to avoid stability issues. 

Our study motivates a challenging question: what does the loss landscape look like? There are some known studies about the loss landscape of a few commonly used loss functions \cite{Loss_Surface}. But to the best of our knowledge, the loss landscape of a loss function that involves the norm of a differential operator of the neural network has not been investigated. If the neural network can well approximate PDE solutions in $H^1$ norm as suggested by \cite{PDE_Approx_NN, weinan2022some}, the "bottom" of a loss function given by the norm of a differential operator should be like a very high dimensional valley, because any boundary condition (resp. initial and boundary condition) can uniquely decide the solution of an elliptic (resp. parabolic) PDE. The neural network training process aims to find a local minimum in this "valley" that also matches the initial distribution, the boundary condition, or the Monte Carlo approximation in our paper. However, the gradient of the loss function orthogonal to those "valleys" may have qualitative difference between a loss function given by an elliptic operator and a loss function given by a parabolic operator. This is because the second order derivative of the neural network approximation is likely being much more sensitive against a random change of connection weights than the first order derivatives. Lack of second order derivative in some directions makes the "valley" less steep. We believe this could be the root cause that makes the Alternating Adam method less effective for time dependent Fokker-Planck equations. We also find that the Alternating Adam method fails frequently when a stationary Fokker-Planck equation has degenerate elliptic term. This further supports our conjecture.

Currently it appears that the best training method for a PINN-like problem is very problem specific.  Alternating Adam works the best for stationary Fokker-Planck equation. The Gradient-Based Momentum Weight method works well for the time dependent Fokker-Planck equation. Many PINN are trained by second order methods such as BFGS \cite{BFGS}. The Trainable Weight method works well for some applications of PINN \cite{Trainable_Weight, Trainable_Weight_2}. However, there is no theory that supports the selection of training methods. After writing this paper, we believe the choice of suitable training method should be dependent on properties of the loss surface. We will address this in our future work.

\section{Author's Contributions, Acknowledgements, and Data Availability Statement}

\subsection{Author's Contributions}

{\bf YL}: Conceptualization, Methodology, Writing - Review and Editing. {\bf CM}: Programming, Analysis, Visualization, Writing - Original Draft, Writing - Review and Editing

\subsection{Acknowledgements} YL is supported by NSF grants DMS-1813246 and DMS-2108628. CM is supported by the REU part of NSF DMS-1813246.

We thank Prof. George Karniadakis and Dr. Shengze Cai for helpful discussions about neural network training, particularly the use of trainable weight.  

\subsection{Data Availability Statement}
The data that support the findings of this study are available from the corresponding author upon reasonable request.

\appendix

\section{Algorithm implementation, hyper-parameter and training point selection}

The majority of this appendix covers the implementation details, hyper-parameter selection, and training data sampling based on the 1D SDE given by equation \eqref{1D} with the multimodal initial distribution. For each training, the error is calculated in the $L^2$ sense in comparison to the ground truth, which is a numerical solution obtained using a Crank-Nicholson solver. The Standard Method is used for selecting training points.  The final subsection covers training point count selection for the Anchor Sampling method. 
 
In this example we used 1500 collocation points with a breakdown of $N^X=500$, $N^Y=1000$, $I^Y=500$. Additionally, $\bm{\mathfrak{Y}}$ is composed of 20,000 points uniformly sampled from the entire domain.  The Monte Carlo sampling part runs $10^7$ samples to approximate the probability density at non-initial collocation points. 

\subsection{Alternating Adam}
While Alternating Adam doesn't have any hyper-parameters to select, a slight modification must be made based on the performance for this SDE.  Unlike for the stationary case, the loss of a given epoch is not a good indicator of performance. The lowest loss is usually seen in the first few epochs. Because of this, the last epoch of training rather that the one with the minimum loss value is selected as the solution of the neural network solver. 

\subsection{Fixed Weight}

The only hyper-parameter that must be selected here is the weight term $\theta$, which serves as a benchmark when comparing the results of non Fixed Weight training methods. Generally $L_1^{\textrm{loss}}$ is much larger than $L_2^{\textrm{loss}}$, so $\theta$ must be close to 1. In the first numerical example we tested the performance of the neural network solver for $10$ different values of $\theta$ in the interval $[0.95, 0.995]$. This interval was selected after a few preliminary tests, which showed that the training result outside of this interval is less satisfactory in general. 
After training with each $\theta$ value in this interval 5 times, the optimal choice was found to be $\theta = 0.975$. 

\begin{figure}[h]
\label{fixed}
{\includegraphics[width=\textwidth]{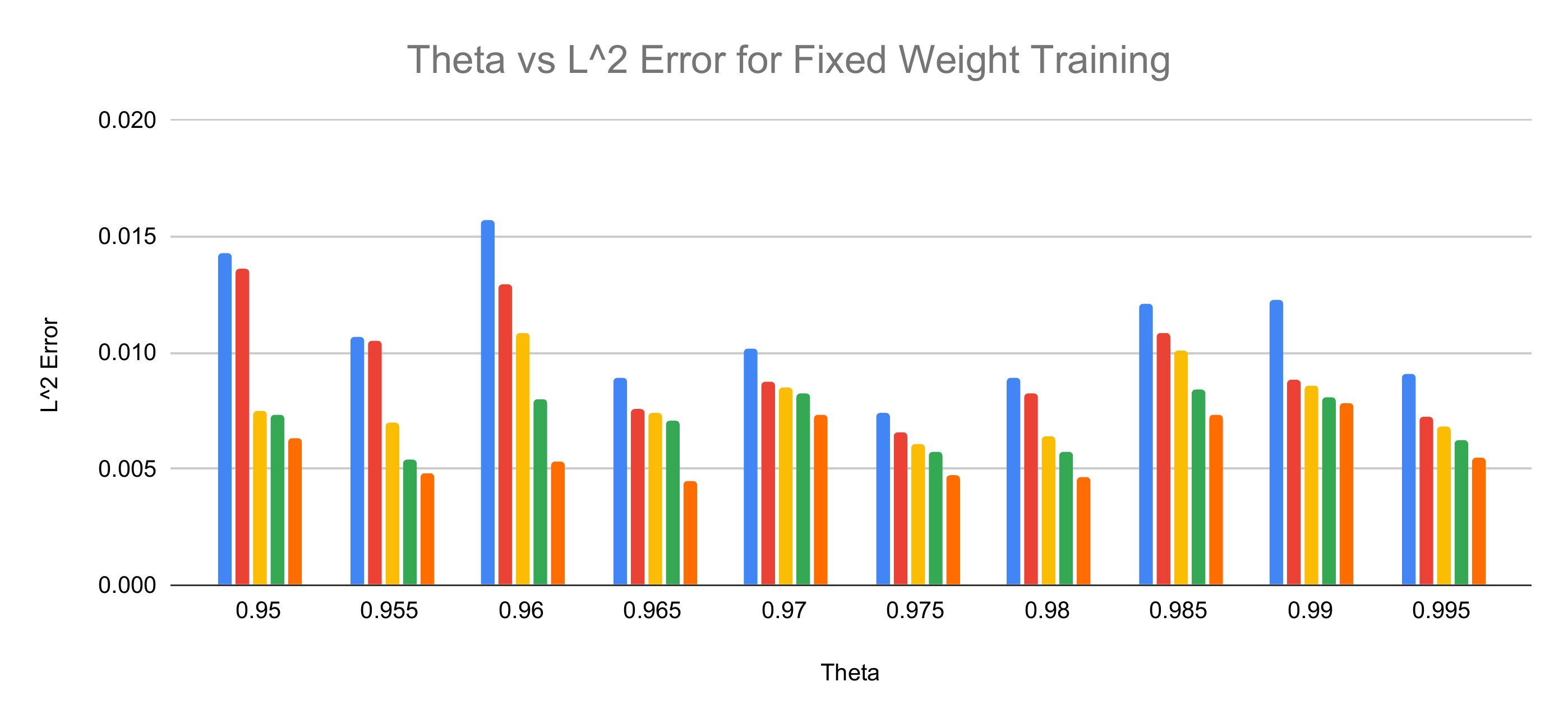}}
\caption{$L^2$ error for different $\theta$ values using Fixed Weight training.}
\end{figure}

This weight is further confirmed by evaluating $45$ weight ratios $\frac{L_1^{\textrm{loss}}}{L_1^{\textrm{loss}}+L_2^{\textrm{loss}}}$ using $\theta = 0.975$. One can see that the $L^2$ error is the lowest when the final weight ratio is between 0.96 and 0.97, and from the results from Figure \ref{fixed} we observed that the final loss ratio averaged around 0.012 less than the $\theta$ value used.

\begin{figure}[t]
\label{006}
{\includegraphics[width=\textwidth]{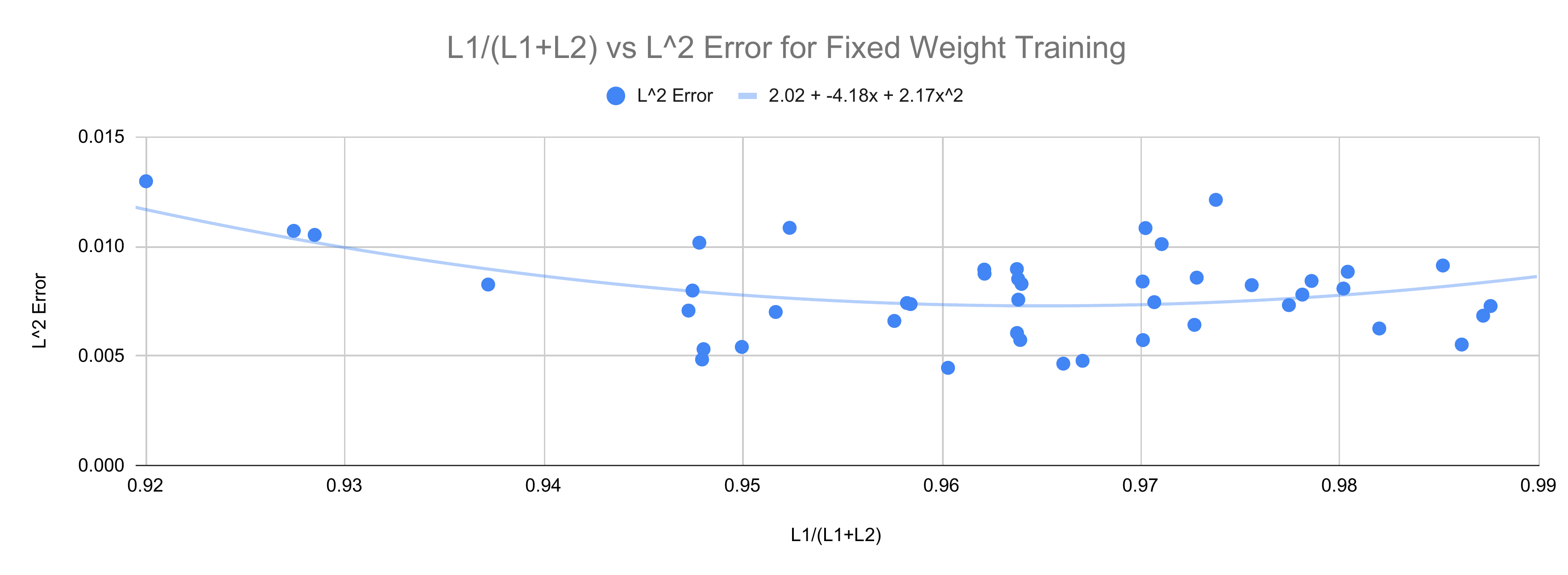}}
\caption{Comparison of final $\frac{L_1^{\textrm{loss}}}{L_1^{\textrm{loss}}+L_2^{\textrm{loss}}}$ and $L^2$ error for all 45 Fixed Weight training results using $\theta = 0.975$.}
\end{figure}

\subsection{Trainable Weight}

The Trainable Weight algorithm only requires selecting an initial value for $\theta_0$.  The simplest choice is $\theta_0=0$. The training results from this initial value can be seen below in Figure \ref{fig2}.  The left panel shows the value of $\theta$ versus the epoch, and the right panel shows the equivalent $\theta$ values if translated to the Fixed Weight algorithm.  As $\theta$ increases roughly linearly, the Fixed Weight equivalent $\theta$ asymptotically approaches 1.  Recall that previously we found that the optimal Fixed Weight $\theta$ is 0.975. This means that the Trainable Weight algorithm quickly moves $\theta$ into the neighborhood of this optimal value then passes this optimal value. Although theoretically the minimax weighting seeks to find a saddle point in the weight space \cite{Trainable_Weight_Theory}, throughout our study, we have not seen the stabilization of $\theta$ as theoretically predicted. As discussed in \cite{Trainable_Weight}, the stabilization may be related to the use of a suitable "mask function". 

\begin{figure}[h]
\label{fig2}
{\includegraphics[width=\textwidth]{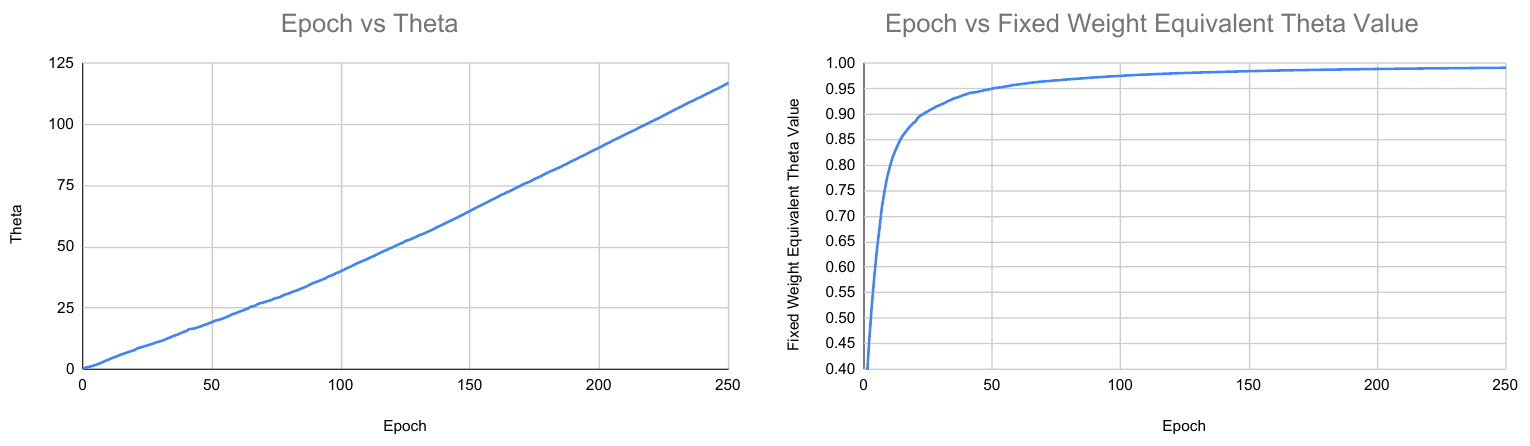}}

\caption{Left: Theta per Epoch for Trainable Weight using $\theta_0=0$.  Right: Equivalent Theta for Fixed Weight per Epoch.}
\end{figure}



\subsection{Loss-Based Momentum Weight} 
\begin{figure}[h]
\label{004}
{\includegraphics[width=\textwidth]{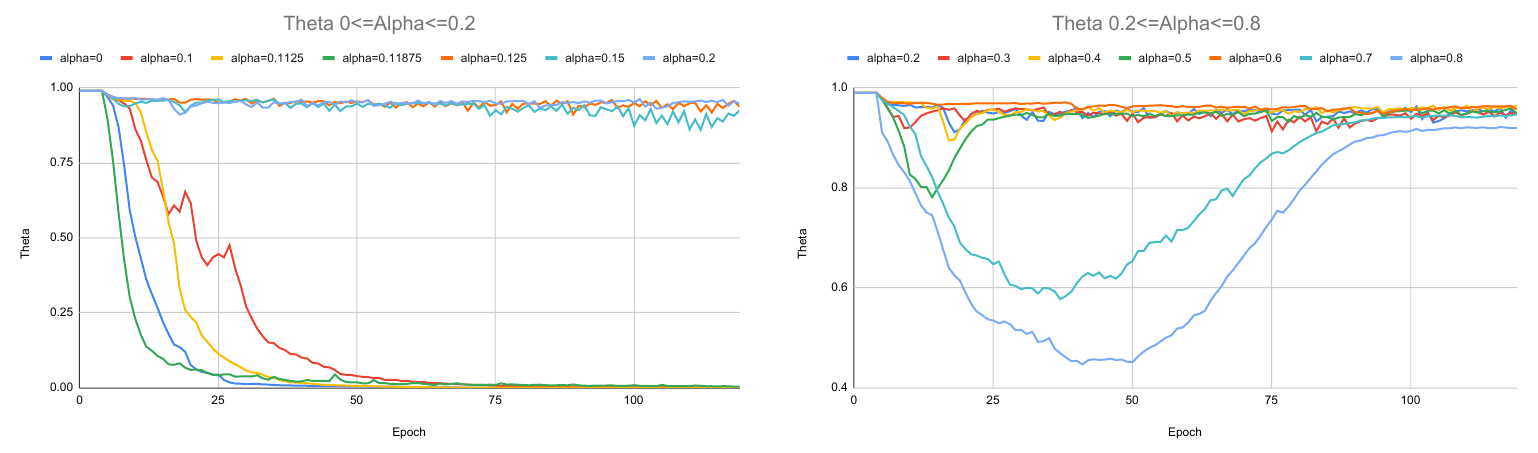}}
\caption{Loss-based Momentum Weight $\theta$ value per epoch for different $\alpha$ values using Loss-Based Momentum Weight. Left: $0 \leq \alpha \leq 0.2$.  Right: $0.2 \leq \alpha \leq 0.8$}
\end{figure}

The Loss-Based Momentum Weight algorithm requires us to choose both the initial value $\theta_0$ and the weight $\alpha$. Since $L_1^{\mathrm{loss}}$ is usually a large quantity in the early phase of training, we chose a $\theta_0$ that is above the expected range of optimal $\theta$ values. Based on some preliminary testing, $\theta_0=0.99$ was selected. The weight $\alpha$ also significantly determines the training dynamics. In our numerical tests the initial loss ratio comes from the 5th epoch (instead of the average of the first $5$ epochs). Figure \ref{004} shows how $\theta$ changes over the epochs depending on the different $\alpha$ values.  This is split into two graphs, since there are two different phenomena present at the low and high range of the $\alpha$ values. 

 The left panel of Figure \ref{004} shows that small values of $\alpha$ make $\theta$ converge towards zero quickly. This is expected because the value of $\theta$ is not properly stabilized, hence the problem of the negative feedback loop is not prevented. The transition starts when $\alpha$ reaches approximately $0.125$, at which point the dynamics of $\theta$ stabilize. On the right panel of Figure \ref{004}, we see a very different behavior going on for larger values of $\alpha$. When $\alpha$ is greater than $0.5$, it becomes harder to revert the initial decreasing trend of $\theta$ because the updates on $\theta$ are too small. As seen in the figure, $\theta$ takes large dips that are followed by increases back towards the stable area near $\theta = 0.95$. Larger values of $\alpha$ make the stabilization even slower. In summary, we conclude that the optimal range of $\alpha$ values appears to be roughly between 0.2 and 0.4.  As larger values make convergence of $\theta$ to 0 less likely, $\alpha=0.4$ was selected.

\begin{figure}[h]
\label{hist}
{\includegraphics[width=\textwidth]{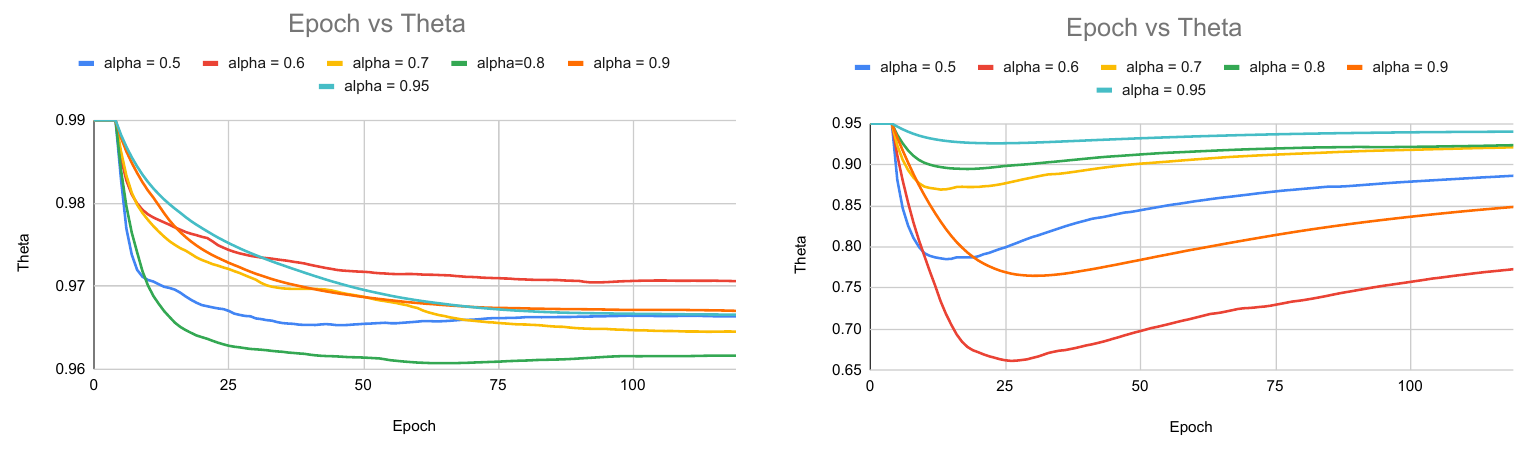}}
\caption{Loss-based Momentum Weight Alternative implementation $\theta$ per epoch for different $\theta_0$ and $\alpha$ values.}
\end{figure}

 As mentioned previously, the alternative implementation of Loss-Based Momentum Weight method uses the average of all historical loss ratios, hence both $\theta_0$ and $\alpha$ needed to be selected again. Figure \ref{hist} demonstrates the evolution of $\theta$ across the epochs when starting from $\theta_0 = 0.99$ (left panel) and $\theta_0 = 0.95$ (right panel). As seen in Figure \ref{hist}, the dynamics of $\theta$ are relatively stable for a wide range of values of $\alpha$, although $\theta$ moves too slowly when $\alpha$ is close to $1$. Here we believe $\theta_0 = 0.99$ remains a good choice of the initial value, and the $L^2$ error results indicate that $\alpha=0.6$ is optimal.

 \subsubsection{Gradient-Based Momentum Weight}
\begin{figure}[h!]
 \label{007}
{\includegraphics[width=\textwidth]{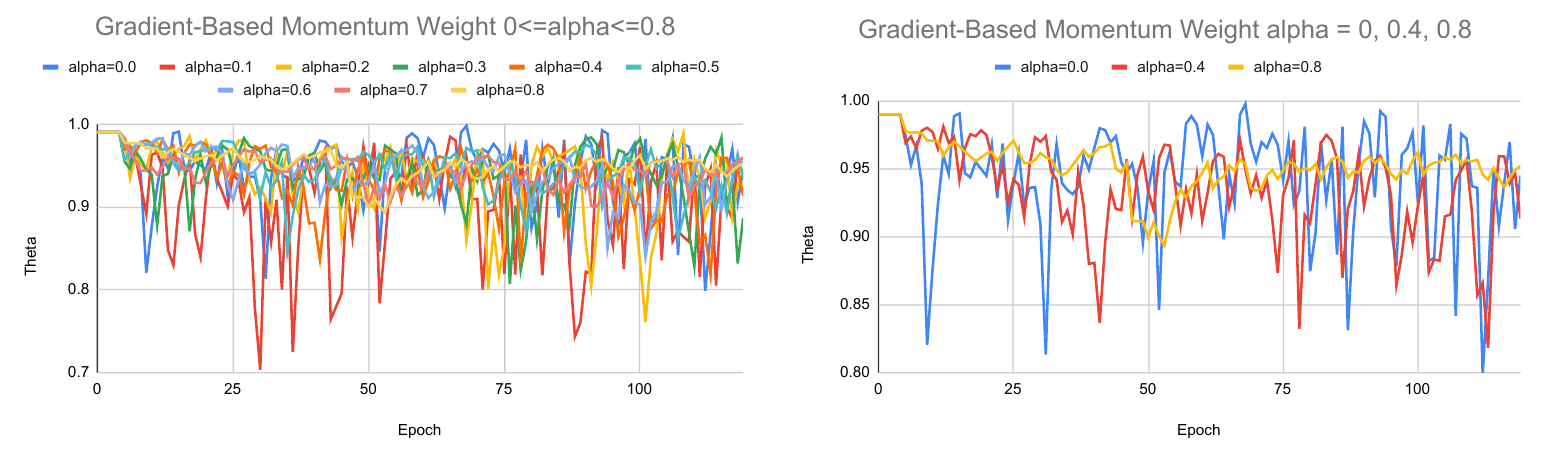}}
\caption{Gradient-Based Momentum Weight $\theta$ value per epoch for Gradient-Based Momentum Weight using different $\alpha$ values.}
\end{figure}

 Finally, Gradient-Based Momentum Weight also requires selecting $\theta_0$ and $\alpha$. Our numerical experiment shows that the performance is not very sensitive against the choice of these hyper-parameters, possibly because the dynamics of $\theta$ have significant fluctuation anyway.   For consistency $\theta_0=0.99$ was selected.  Figure \ref{007} shows $\theta$ per epoch for different $\alpha$ values.  Increasing $\alpha$ decreases the variance and vice versa, but no other behavior is introduced by changing $\alpha$, and performance remains similar across all $\alpha$ values. Further tests showed that $L^2$ error is lowest at $\alpha=0.4$, so this was selected. 
 
\subsection{Anchor Sampling Implementation}

The median and mean $L^2$ error over the time domain for five training point counts at $t=0.2$ can be seen below in Figure $\ref{anchor}$.  Grid selection and 120 epochs of training were used here.  1156 and 2500 points produced similar results, while the rest produced worse results.  Because of this 1156 was selected.   

\begin{figure}[h!]
 \label{anchor}
{\includegraphics[width=\textwidth]{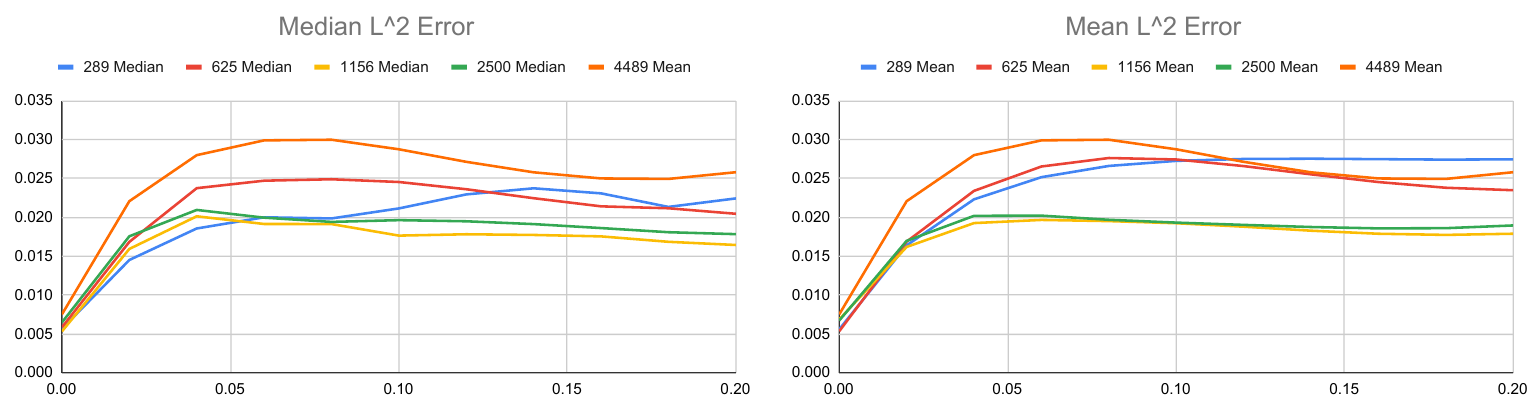}}
\caption{Median and Mean $L^2$ error from 11 trainings using Anchor Sampling with Grid Selection for 120 epochs.  Sample counts of 289, 625, 1156, 2500, and 4489 are used at $t=0.2$.}
\end{figure}

\bibliographystyle{plain}
\bibliography{FPEref.bib}

\begin{thebibliography}{10}

\bibitem{FPEPINN}
Xiaoli Chen, Liu Yang, Jinqiao Duan, and George~Em Karniadakis.
\newblock Solving inverse stochastic problems from discrete particle
  observations using the fokker--planck equation and physics-informed neural
  networks.
\newblock {\em SIAM Journal on Scientific Computing}, 43(3):B811--B830, 2021.

\bibitem{PDE_Approx_NN}
Ziang Chen, Jianfeng Lu, and Yulong Lu.
\newblock On the representation of solutions to elliptic pdes in barron spaces.
\newblock In Marc'Aurelio Ranzato, Alina Beygelzimer, Yann~N. Dauphin, Percy
  Liang, and Jennifer~Wortman Vaughan, editors, {\em Advances in Neural
  Information Processing Systems 34: Annual Conference on Neural Information
  Processing Systems 2021, NeurIPS 2021, December 6-14, 2021, virtual}, pages
  6454--6465, 2021.

\bibitem{FPE_2022}
Matthew Dobson, Yao Li, and Jiayu Zhai.
\newblock An efficient data-driven solver for fokker–planck equations:
  Algorithm and analysis.
\newblock {\em Communications in Mathematical Sciences}, 20(3):803–827, 2022.

\bibitem{Noisy_PINN}
Hamidreza Eivazi and Ricardo Vinuesa.
\newblock Physics-informed deep-learning applications to experimental fluid
  mechanics.
\newblock {\em arXiv preprint arXiv:2203.15402}, 2022.

\bibitem{Trainable_Weight}
Yiqi Gu, Haizhao Yang, and Chao Zhou.
\newblock Selectnet: Self-paced learning for high-dimensional partial
  differential equations.
\newblock {\em Journal of Computational Physics}, 441:110444, 2021.

\bibitem{Adam}
Diederik~P Kingma and Jimmy Ba.
\newblock Adam: A method for stochastic optimization.
\newblock {\em arXiv preprint arXiv:1412.6980}, 2014.

\bibitem{Loss_Surface}
Daniel~D Lee, P~Pham, Y~Largman, and A~Ng.
\newblock Advances in neural information processing systems 22.
\newblock Technical report, Tech. Rep., Tech. Rep, 2009.

\bibitem{FPE_2019}
Yao Li.
\newblock A data-driven method for the steady state of randomly perturbed
  dynamics.
\newblock {\em Communications in Mathematical Sciences}, 17(4):1045–1059,
  2019.

\bibitem{Trainable_Weight_Theory}
Dehao Liu and Yan Wang.
\newblock A dual-dimer method for training physics-constrained neural networks
  with minimax architecture.
\newblock {\em Neural Networks}, 136:112--125, 2021.

\bibitem{BFGS}
Dong~C Liu and Jorge Nocedal.
\newblock On the limited memory bfgs method for large scale optimization.
\newblock {\em Mathematical programming}, 45(1):503--528, 1989.

\bibitem{Trainable_Weight_2}
Levi McClenny and Ulisses Braga-Neto.
\newblock Self-adaptive physics-informed neural networks using a soft attention
  mechanism.
\newblock {\em arXiv preprint arXiv:2009.04544}, 2020.

\bibitem{Karniadakis_PINN}
Maziar Raissi, Paris Perdikaris, and George~E Karniadakis.
\newblock Physics-informed neural networks: A deep learning framework for
  solving forward and inverse problems involving nonlinear partial differential
  equations.
\newblock {\em Journal of Computational physics}, 378:686--707, 2019.

\bibitem{thomas2013numerical}
James~William Thomas.
\newblock {\em Numerical partial differential equations: finite difference
  methods}, volume~22.
\newblock Springer Science \& Business Media, 2013.

\bibitem{weinan2022some}
E~Weinan and Stephan Wojtowytsch.
\newblock Some observations on high-dimensional partial differential equations
  with barron data.
\newblock In {\em Mathematical and Scientific Machine Learning}, pages
  253--269. PMLR, 2022.

\bibitem{NN_FPE}
Jiayu Zhai, Matthew Dobson, and Yao Li.
\newblock A deep learning method for solving fokker-planck equations.
\newblock In Joan Bruna, Jan Hesthaven, and Lenka Zdeborova, editors, {\em
  Proceedings of the 2nd Mathematical and Scientific Machine Learning
  Conference}, volume 145 of {\em Proceedings of Machine Learning Research},
  pages 568--597. PMLR, 16--19 Aug 2022.

\end{thebibliography}
\end{document}